\documentclass{amsproc}%
\usepackage{amsmath}
\usepackage{amsfonts}
\usepackage{amssymb}
\usepackage{graphicx}%
\newtheorem{theorem}{Theorem}[section]

\theoremstyle{definition}

\newtheorem{notation}[theorem]{Notation}
\theoremstyle{remark}

\newtheorem{lem}[theorem]{Lemma}
\newtheorem{cor}[theorem]{Corollary}
\newtheorem{exmp}[theorem]{Example}
\newtheorem{prop}[theorem]{Proposition}
\newtheorem{rem}[theorem]{Remark}
\newtheorem{thm}[theorem]{Theorem}
\newenvironment{pf}{\begin{proof}}{\end{proof}}
\numberwithin{equation}{section}
\begin{document}
\title{Certain representations of the Cuntz relations, and a question on wavelets decompositions}
\author{Palle E.T. Jorgensen}
\address{Department of Mathematics, The University of Iowa, 14 MacLean Hall, Iowa City,
IA, 52242-1419, U.S.A.}
\email{http://www.math.uiowa.edu/\symbol{126}jorgen}
\thanks{This material is based upon work supported by the U.S. National Science
Foundation under Grant No. DMS-0139473 (FRG)}
\date{}

\begin{abstract}
We compute the Coifman-Meyer-Wickerhauser measure $\mu$ for certain families
of quadrature mirror filters (QMFs), and we establish that for a subclass of
QMFs, $\mu$ contains a fractal scale.

\end{abstract}
\maketitle

\section{Introduction\label{Intro}}

\markboth{PALLE E. T. JORGENSEN}{REPRESENTATIONS OF CUNTZ REL. / WAVELETS DECOMPS.}
It is known that certain representations of the Cuntz relations (or
equivalently the Cuntz algebras $\mathcal{O}_{N}$ \cite{Cu77} serve as a
computational tool for wavelet analysis. The reason for this is that a
representation of the Cuntz relations on a Hilbert space $\mathcal{H}$ induces
a successive subdivision of $\mathcal{H}$ into orthogonal subspaces. To get
the scaling feature of wavelets enter into our representations, build a family
of operators from the endomorphism $z\longrightarrow z^{N}$ on the circle
$\mathbb{T}$, and $N$ carefully chosen multiplication operators on
$L^{2}(\mathbb{T})$ with respect to Haar measure on $\mathbb{T}$. In the
dyadic case, we work with $z\longrightarrow z^{2}$ and with two multiplication
operators, i.e., two functions on $\mathbb{T}$. In the language of signal
processing, these functions $m_{0}$, and $m_{1}$ are called low-pass and
high-pass filters, respectively; or the pair is called a quadrature-mirror
system of filters. By a standard procedure \cite{Dau92}, we construct the
scaling function (or father function) $\varphi$ on $\mathbb{R}$ from $m_{0}$,
and the wavelet function (or mother function) $\psi$ from $m_{1}$. We get two
corresponding isometries which yield a representation of the Cuntz algebra
$\mathcal{O}_{2}$. These representations have been studied in earlier papers,
see e.g., \cite{BJMP04, BraEvaJor2000, BrJo02b, BrJo99, DuJo05, JoKr03}. The
problem in wavelet theory is to build orthonormal bases in $L^{2}(\mathbb{R})$
from these data. This can be done \cite{Dau92}, and the wavelet bases have
advantages over the earlier known basis constructions, one of the main
advantages is efficiency of computation.\addtolength{\textheight}{-12pt}

While choice pyramids (see e.g., \cite{Wic93} and \cite{Wic94}) have been much
used in computations already, the presence of representations of
$\mathcal{O}_{N}$ yields a Hilbert space geometric way of realize such
combinatorial pyramids. Representing the integers $\mathbb{N}_{0}$
$=\{0,1,2,\ldots\}$ using the Euclidean algorithm, and using a corresponding
pyramid, we get a orthogonal basis for $L^{2}(\mathbb{R}),\varphi_{0}%
,\varphi_{1},\varphi_{2},\ldots$ called a \textit{wavelet packet}. A wavelet
packet decomposition arises by a combinatorial selection of a subset of the
index set $\mathbb{N}_{0}$ for these basis function, and a corresponding
subset of the scaling operations $x\longrightarrow N^{j}x$ (i.e., a selection
of a subset of the integers $\mathbb{Z}$) such that this selection of subsets
yields an orthonormal basis for $L^{2}(\mathbb{R})$. As pointed out in
\cite{Wic93} and \cite{Wic94}, the selection of such feasible pairs of subsets
involves a combinatorial tiling question, and a problem in measure theory. In
this paper we show how these questions may conveniently be addressed with the
use of the non-commutative harmonic analysis coming from the theory of
representations of the $C^{\ast}$-algebra $\mathcal{O}_{N}$.

And in the course of the paper, we show that the scope of our method is wider
than the original context of standard wavelets: It applies equally well to the
class of fractals that arise from affine iterated function systems (IFS); see
e.g., \cite{DuJo05}, \cite{Falconer85}, and \cite{Hut81}.

Wavelets and wavelet packets are orthonormal bases, or Parseval frames, for
the Hilbert space $L^{2}(\mathbb{R})$ which are built from a small set of
generating functions by applying only two operations to them, translation by
the integers $\mathbb{Z\subset R}$, and by dyadic scaling $t\longrightarrow
2^{j}t,j\in\mathbb{Z}$. The best known case is that of singly generated
wavelets. Then there is a function $\psi\in L^{2}(\mathbb{R})$ such that the
family
\begin{equation}
\left\{  2^{j/2}\psi\left(  2^{j}t-k\right)  \right\}  _{j,k\in\mathbb{Z}}
\label{IntroEq1}%
\end{equation}
is an orthonormal basis for $L^{2}\left(  \mathbb{R}\right)  $. Wavelet
packets (in the orthonormal case) are sequences $\left(  \varphi_{n}\right)
_{n\in\mathbb{N}_{0}}$ and subsets $J\subset\mathbb{N}_{0}\times\mathbb{Z}$
such that the two families
\begin{equation}
\left\{  2^{j/2}\varphi_{n}\left(  2^{j}t-k\right)  \mid\left(  n,j\right)
\in J,k\in\mathbb{Z}\right\}  \label{IntroEq2}%
\end{equation}
and
\begin{equation}
\left\{  \varphi_{n}\left(  t-k\right)  \mid n\in\mathbb{N}_{0},k\in
\mathbb{Z}\right\}  \label{IntroEq3}%
\end{equation}
are both both orthonormal bases.

In the seminal paper \cite{CMW92}, the authors Coifman, Meyer, and
Wickerhauser proposed a selection of systems (\ref{IntroEq2}) from
(\ref{IntroEq3}); the idea being that (\ref{IntroEq3}) may be constructed by
an effective and relatively simple algorithm. The construction starts with a
quadrature-mirror wavelet filter, i.e., a Fourier series
\begin{equation}
m_{0}(x)=\sum_{k\in\mathbb{Z}}a_{k}e^{-i2\pi kx} \label{IntroEq4}%
\end{equation}
such that
\begin{equation}
\sum_{k\in\mathbb{Z}}\overline{a_{k}}a_{k+2\ell}=\delta_{0,\ell}
\label{IntroEq5}%
\end{equation}
and%
\begin{equation}
\sum_{k\in\mathbb{Z}}a_{k}=\sqrt{2}\text{.} \label{IntroEq6}%
\end{equation}
If the expansion in (\ref{IntroEq4}) is a finite sum, then the functions
$\varphi_{0},\varphi_{1},\varphi_{2},\cdots$ will be of compact support. The
advantage with using (\ref{IntroEq2}) over (\ref{IntroEq1}) or (\ref{IntroEq3}%
) is that the basis functions in (\ref{IntroEq2}) are better localized in
time-frequency, as is spelled out in \cite{CMW92} and \cite{Wic93}.

The authors of \cite{CMW92} propose a certain integral decomposition of
$L^{2}\left(  \mathbb{R}\right)  $ which is based on a certain measure
$\mu_{0}$ on the unit-interval $\left[  0,1\right]  $. The idea is the
following: The quadrature rules (\ref{IntroEq5}) lead to a dyadic
decomposition of the Hilbert space $\ell^{2}\left(  \mathbb{Z}\right)  $,
repository for the wavelet coefficients: For each $k\in\mathbb{N}$, and each
dyadic interval
\begin{equation}
\left[  \xi,\xi+2^{-k}\right)  \text{, }\xi=\frac{i_{1}}{2}+\cdots+\frac
{i_{k}}{2^{k}}\text{,} \label{IntroEq7}%
\end{equation}
we assign a closed subspace $\mathcal{H}\left(  \xi\right)  \subset\ell^{2}$,
such that the $2^{k}$ distinct subspaces $\mathcal{H}\left(  \xi\right)  $ are
mutually orthogonal, and
\begin{equation}
\sum_{\xi}^{\oplus}\mathcal{H}\left(  \xi\right)  =\ell^{2} \label{IntroEq8}%
\end{equation}
If $e_{0}$ denotes the vector $e_{0}\left(  j\right)  $:$=\delta_{0,j}$ in
$\ell^{2}$, and $P\left(  \xi\right)  $ the orthogonal projection onto
$\mathcal{H}\left(  \xi\right)  $, then the measure $\mu_{0}$ is determined
by
\begin{equation}
\mu_{0}\left(  \left[  \xi,\xi+2^{-k}\right)  \right)  =\left\Vert P\left(
\xi\right)  e_{0}\right\Vert ^{2} \label{IntroEq9}%
\end{equation}
where $\xi$ is the dyadic rational (\ref{IntroEq7}).

More generally, if $f\in\mathcal{H}$, $\left\Vert f\right\Vert ^{2}=1,$ we set
$\mu_{f}\left(  \cdot\right)  =\left\Vert P\left(  \cdot\right)  f\right\Vert
^{2}$, and note that each $\mu_{f}$ is a probability measure. These measures
are shown in section \ref{Measures} to dictate the wavepacket analysis in the
Hilbert space $L^{2}\left(  \mathbb{R}\right)  $.

The expectation is that $\mu_{0}$ is absolutely continuous with respect to the
Lebesgue measure on $\left[  0,1\right]  $. We give a formula for $\mu_{0}$
directly in terms of the coefficients $\left(  a_{j}\right)  _{j\in\mathbb{Z}%
}$ in (\ref{IntroEq4}), and we give a one-parameter family of coefficients
$a_{j}\left(  \beta\right)  $ where $\beta$ parameterizes the circle, and
where for $\left\vert \beta\right\vert <\pi/4$, the measure $\mu_{0}$ on
\textquotedblleft small\textquotedblright\ dyadic intervals $J$ is governed by
the formula
\begin{equation}
\mu_{0}\left(  \text{dyadic interval }J\right)  \approxeq\left\vert
J\right\vert ^{s}\text{,} \label{IntroEq10}%
\end{equation}
where $\left\vert J\right\vert $ denotes the length of $J$, and where
\begin{equation}
s=\frac{\ln\left(  a_{0}\left(  \beta\right)  ^{-2}\right)  }{\ln2}\text{.}
\label{IntroEq11}%
\end{equation}
In the interval for $\beta$, the coefficients $a_{j}\left(  \beta\right)  $
are real-valued, and when $\left\vert \beta\right\vert <\pi/4$, $a_{0}\left(
\beta\right)  >1/\sqrt{2}$.

\section{Fourier Polynomials\label{FourPoly}}

Let $\mathbb{T}=\left\{  z\in\mathbb{C}\mid\left\vert z\right\vert =1\right\}
$, and let $D\in\mathbb{N}$. Consider the following two functions
\begin{equation}
m_{0}\left(  z\right)  =\sum_{k=0}^{2D-1}a_{k}z^{k} \label{FourEq1}%
\end{equation}
and%
\begin{align}
m_{1}\left(  z\right)   &  =z^{2D-1}\overline{m_{0}\left(  -z\right)
}\label{FourEq2}\\
&  =\sum_{k=0}^{2D-1}\left(  -1\right)  ^{k}\bar{a}_{2D-1-k}z^{k}\nonumber
\end{align}
We shall need several lemmas.

\begin{lem}
\label{FourLem1}The finite sequence $a_{0},a_{1},\cdots,a_{2D-1}$ in
\emph{(}\ref{FourEq1}\emph{) }satisfies
\begin{equation}
\sum_{j\in\mathbb{Z}}\bar{a}_{k}a_{k+2\ell}=\delta_{0,\ell} \label{FourEq3}%
\end{equation}
if and only if anyone of the following equivalent conditions holds:

\noindent\textup{(}a\textup{)} The matrix function
\begin{equation}
U\left(  z\right)  =\sum_{k=0}^{D-1}A_{k}z^{k}\text{, }z\in\mathbb{T}
\label{FourEq4}%
\end{equation}
is unitary, where the coefficient matrices are
\begin{equation}
A_{k}=\left(
\begin{array}
[c]{cc}%
a_{2k} & a_{2k+1}\\
a_{2\left(  D-k\right)  -1} & -a_{2\left(  D-k-1\right)  }%
\end{array}
\right)  \text{ for }k=0,1,\cdots,D-1\text{.} \label{FourEq5}%
\end{equation}
\noindent\textup{(}b\textup{)} The operators
\begin{equation}
\left(  S_{i}f\right)  \left(  z\right)  =m_{i}\left(  z\right)  f\left(
z^{2}\right)  \text{, }f\in L^{2}\left(  \mathbb{T}\right)  \text{,
}i=0,1\text{,} \label{FourEq6}%
\end{equation}
and their adjoints $S_{i}^{\ast}$ satisfy
\begin{equation}
\left\{
\begin{array}
[c]{c}%
S_{i}^{\ast}S_{j}=\delta_{i,j}I\text{, }i,j=0,1\\
S_{0}S_{0}^{\ast}+S_{1}S_{1}^{\ast}=I
\end{array}
\text{.}\right.  \label{FourEq7}%
\end{equation}
\noindent\textup{(}c\textup{)} The following matrix function
\begin{equation}
M\left(  z\right)  \text{\emph{:}}=\frac{1}{\sqrt{2}}\left(
\begin{array}
[c]{cc}%
m_{0}\left(  z\right)  & m_{0}\left(  -z\right) \\
m_{1}\left(  z\right)  & m_{1}\left(  -z\right)
\end{array}
\right)  \text{, }z\in\mathbb{T}\text{,} \label{FourEq8}%
\end{equation}
is unitary, i.e., $M\left(  z\right)  ^{\ast}M\left(  z\right)  =I$, for all
$z\in\mathbb{T}$.
\end{lem}

\begin{rem}
\label{FourRem1}In the summation \textup{(}\ref{FourEq3}\textup{)}, it is
understood that the terms are zero if the summation index is not in the range
which is specified for the sequence: If for example $D=2$, then the conditions
\textup{(}\ref{FourEq3}\textup{)} spell out as follows:
\begin{equation}
\left\{
\begin{array}
[c]{c}%
\left\vert a_{0}\right\vert ^{2}+\left\vert a_{1}\right\vert ^{2}+\left\vert
a_{2}\right\vert ^{2}+\left\vert a_{3}\right\vert ^{2}=1\\
\bar{a}_{0}a_{2}+\bar{a}_{1}a_{3}=0
\end{array}
\text{.}\right.  \label{FourEq9}%
\end{equation}
The interpretation of \textup{(}\ref{FourEq4}\textup{)} is that the matrices
$A_{0},A_{1},\cdots A_{D-1}$ are Fourier coefficients for a matrix function
$z\longrightarrow U\left(  z\right)  $ defined on $\mathbb{T}$. The
requirements is that for each $z\in\mathbb{T}$, the 2 by 2 matrix $U\left(
z\right)  $ is unitary, i.e., that
\begin{equation}
U(z)^{\ast}U\left(  z\right)  =I\text{, }z\in\mathbb{T}\text{.}
\label{FourEq10}%
\end{equation}
In general, for a function of the form \textup{(}\ref{FourEq4}\textup{)}, the
unitarity condition \textup{(}\ref{FourEq10}\textup{)} is equivalent to the
following conditions on the constant matrices $A_{0},A_{1},\cdots$:
\begin{equation}
\sum_{k}A_{k+n}A_{k}^{\ast}=\delta_{0,n} \label{FourEq11}%
\end{equation}
or equivalently
\begin{equation}
\sum_{k}A_{k}^{\ast}A_{k+n}=\delta_{0,n}\text{.} \label{FourEq12}%
\end{equation}
Again if $D=2$, these conditions flesh out as follows: The two matrices
\begin{equation}
A_{0}=\left(
\begin{array}
[c]{cc}%
a_{0} & a_{1}\\
a_{3} & -a_{2}%
\end{array}
\right)  \text{ and }A_{1}=\left(
\begin{array}
[c]{cc}%
a_{2} & a_{3}\\
a_{1} & -a_{0}%
\end{array}
\right)  \label{FourEq13}%
\end{equation}
satisfy
\begin{equation}
\left\{
\begin{array}
[c]{c}%
A_{0}A_{0}^{\ast}+A_{1}A_{1}^{\ast}=I\\
A_{1}A_{0}^{\ast}=0
\end{array}
\text{.}\right.  \label{FourEq14}%
\end{equation}

\end{rem}

\begin{pf}
\textbf{of Lemma \ref{FourLem1}.} The coefficients $a_{0},a_{1},\cdots
,a_{2D-1}$ are given and the two functions $m_{0}(z)$ and $m_{1}(z)$ are
defined as in (\ref{FourEq1}) and (\ref{FourEq2}). Even though they are
defined for all $z\in\mathbb{C}$, we shall use only the restrictions to
$\mathbb{T}$. It is convenient in some later calculations to introduce the
substitution $z$:$=e^{-i2\pi x}$; and with this substitution we shall also
view $m_{0}$ and $m_{1}$ as one-periodic functions on the real line
$\mathbb{R}$, and make the identification $\mathbb{T}\cong\mathbb{R}%
/\mathbb{Z}$. We now turn to the implications: \noindent\textup{(}%
\ref{FourEq3}\textup{)}$\Longrightarrow$\textup{(}a\textup{)}. We already
noted that the unitarity of the function $U(z)$ in \textup{(}\ref{FourEq4}%
\textup{)} may be expressed by the conditions \textup{(}\ref{FourEq11}%
\textup{)}. When \textup{(}\ref{FourEq5}\textup{)} is substituted in
\textup{(}\ref{FourEq11}\textup{)}, it is immediate that the two identities
\textup{(}\ref{FourEq3}\textup{)} and \textup{(}\ref{FourEq11}\textup{)} are
equivalent. \noindent\textup{(}a\textup{)}$\Longrightarrow$\textup{(}%
b\textup{)}. Using the normalized Haar measure $\mu$ on $\mathbb{T}$ and the
corresponding inner product
\begin{equation}
\left\langle f\mid g\right\rangle \text{:}=\int_{\mathbb{T}}\overline{f\left(
z\right)  }g\left(  z\right)  \;d\mu\left(  z\right)  \label{FourEq15}%
\end{equation}
of the Hilbert space $L^{2}\left(  \mathbb{T}\right)  $, we get the formula
\begin{equation}
\left(  S_{i}^{\ast}f\right)  \left(  z\right)  =\frac{1}{2}\sum
_{w\in\mathbb{T},w^{2}=z}\overline{m_{i}\left(  w\right)  }f\left(  w\right)
\text{, defined for }f\in L^{2}\left(  \mathbb{T}\right)  \text{ },\text{
}z\in\mathbb{T}\text{ , and }i=0,1\text{.} \label{FourEq16}%
\end{equation}
As a result, we get the formula
\begin{equation}
\left(  S_{i}^{\ast}S_{j}f\right)  \left(  z\right)  =\frac{1}{2}\sum
_{w\in\mathbb{T},w^{2}=z}\overline{m_{i}\left(  w\right)  }m_{j}\left(
w\right)  f\left(  z\right)  \text{;} \label{FourEq17}%
\end{equation}
which states that the four operators $S_{i}^{\ast}S_{j}$ are multiplication
operators. Hence the first part of \textup{(}\ref{FourEq7}\textup{)} reads:
\begin{equation}
\frac{1}{2}\sum_{w\in\mathbb{T},w^{2}=z}\overline{m_{i}\left(  w\right)
}m_{j}\left(  w\right)  =\delta_{i,j}\text{ , }z\in\mathbb{T}\text{.}
\label{FourEq18}%
\end{equation}
Substitution of \textup{(}\ref{FourEq1}\textup{)} into \textup{(}%
\ref{FourEq18}\textup{)} shows that these conditions amount to the same sum
rules \textup{(}\ref{FourEq3}\textup{)}, or equivalently \textup{(}%
\ref{FourEq11}\textup{)}. The second equation from \textup{(}\ref{FourEq7}%
\textup{)} may be restated as
\begin{equation}
\left\Vert S_{0}^{\ast}f\right\Vert ^{2}+\left\Vert S_{1}^{\ast}f\right\Vert
^{2}=\left\Vert f\right\Vert ^{2}\text{ , for all }f\in L^{2}\left(
\mathbb{T}\right)  \text{;} \label{FourEq19}%
\end{equation}
or equivalently
\begin{equation}
\sum_{i=0}^{1}\left\vert \frac{1}{2}\sum_{w^{2}=z}\overline{m_{i}\left(
w\right)  }f\left(  w\right)  \right\vert ^{2}d\mu\left(  z\right)
=\int_{\mathbb{T}}\left\vert f\left(  z\right)  \right\vert ^{2}\text{ }%
d\mu\left(  z\right)  \text{.} \label{FourEq20}%
\end{equation}
But the summation on the left-hand sides fleshes out as follows:
\begin{align*}
&  \frac{1}{4}\sum_{i=0}^{1}\sum_{w^{2}=\zeta^{2}=z}\overline{m_{i}\left(
w\right)  }m_{i}\left(  \zeta\right)  f\left(  w\right)  f\overline{\left(
\zeta\right)  }\\
&  =\frac{1}{2}\sum_{w^{2}=\zeta^{2}=z}\left(  \frac{1}{2}\sum_{i=0}%
^{1}\overline{m_{i}\left(  w\right)  }m_{i}\left(  \zeta\right)  \right)
f\left(  w\right)  f\overline{\left(  \zeta\right)  }\\
&  =\frac{1}{2}\sum_{w^{2}=\zeta^{2}=z}\delta_{w,\zeta}f\left(  w\right)
f\overline{\left(  \zeta\right)  }\\
&  =\frac{1}{2}\sum_{w^{2}=z}\left\vert f\left(  w\right)  \right\vert
^{2}\text{,}%
\end{align*}
and
\begin{equation}
\int_{\mathbb{T}}\frac{1}{2}\sum_{w^{2}=z}\left\vert f\left(  w\right)
\right\vert ^{2}\text{ }d\mu\left(  z\right)  =\int_{\mathbb{T}}\left\vert
f\left(  z\right)  \right\vert ^{2}\text{ }d\mu\left(  z\right)
\label{FourEq21}%
\end{equation}
by an elementary property of the Haar measure. This proves \textup{(}%
a\textup{)}$\Longrightarrow$\textup{(}b\textup{)}; and in fact equivalence.
\noindent\textup{(}b\textup{)}$\Longrightarrow$\textup{(}c\textup{)}. This
implication also goes in both directions, and it is implicit in the previous
step where we showed that the identities \textup{(}\ref{FourEq18}\textup{)}
are equivalent to \textup{(}\ref{FourEq3}\textup{)}. But an inspection of
matrix entries shows that \textup{(}\ref{FourEq18}\textup{)} is a restatement
of the unitary property for the matrix $M\left(  z\right)  $ in \textup{(}%
c\textup{)}.
\end{pf}

\begin{cor}
\label{FourCor1} Consider a Fourier polynomial
\begin{equation}
m\left(  z\right)  =\sum_{k=0}^{2D-1}a_{k}z^{k} \label{FourEq22}%
\end{equation}
as in Lemma \ref{FourLem1}, and let $S=S_{m}$ be the operator on $L^{2}\left(
\mathbb{T}\right)  $ defined by
\begin{equation}
\left(  Sf\right)  \left(  z\right)  =m\left(  z\right)  f\left(
z^{2}\right)  \text{ , }z\in\mathbb{T}\text{.} \label{FourEq23}%
\end{equation}
Then the following conditions are equivalent:

\noindent\textup{(}i\textup{)} $S$ is an isometry;

\noindent\textup{(}ii\textup{)} $\left\vert m\left(  z\right)  \right\vert
^{2}+\left\vert m\left(  -z\right)  \right\vert ^{2}=2$ , for $z\in\mathbb{T}%
$; and

\noindent\textup{(}iii\textup{)} $\sum_{k}\overline{a}_{k}a_{k+2\ell}%
=\delta_{0,\ell}$.
\end{cor}

\begin{pf}
The equivalence \textup{(}i\textup{)}$\Longleftrightarrow$\textup{(}%
ii\textup{)} is immediate from an inspection of \textup{(}\ref{FourEq18}%
\textup{)} in the special case $i=j=0$. To see \textup{(}ii\textup{)}%
$\Longleftrightarrow$\textup{(}iii\textup{)}, it is simplest to write out the
Fourier coefficients of the function
\begin{equation}
z\longrightarrow\frac{1}{2}\sum_{w\in\mathbb{T},\text{ }w^{2}=z}\left\vert
m\left(  w\right)  \right\vert ^{2}\text{.} \label{FourEq24}%
\end{equation}
A substitution of \textup{(}\ref{FourEq22}\textup{)} into \textup{(}%
\ref{FourEq24}\textup{)} yields
\begin{align*}
&  \frac{1}{2}\sum_{j}\sum_{k}\sum_{w^{2}=z}\overline{a_{j}}a_{k}w^{k-j}\\
&  =\sum_{n}\sum_{j,k\,s.t.\,k-j=2n}\overline{a}_{j}a_{k}\left(  \frac{1}%
{2}\sum_{w^{2}=z}w^{2n}\right) \\
&  =\sum_{n}\left(  \sum_{j}\overline{a}_{j}a_{j+2n}\right)  z^{n}\text{,}%
\end{align*}
and it follows that the sum in \textup{(}\ref{FourEq24}\textup{)} is the
constant function $1$ if and only if \textup{(}iii\textup{)} is satisfied.
Note that \textup{(}iii\textup{)} is a restatement of condition \textup{(}%
\ref{FourEq3}\textup{)} in Lemma \ref{FourLem1}. In the following we consider
operators $S$ on the form $(Sf)(z)=m\left(  z\right)  f\left(  z^{2}\right)
$, $z\in\mathbb{T}$, $f\in L^{2}\left(  \mathbb{T}\right)  $, and their
adjoints. But we will need these operators realized explicitly in the sequence
space $\ell^{2}\cong L^{2}\left(  \mathbb{T}\right)  $.
\end{pf}

\begin{lem}
\label{FourLem2}.~

\noindent\textup{(}a\textup{)} The matrix representation for the operator $S$
in Corollary \ref{FourCor1} relative to the standard basis $e_{n}\left(
z\right)  $:$=z^{n}$ , $n\in\mathbb{Z}$ , for $L^{2}\left(  \mathbb{T}\right)
$, is
\begin{equation}
\left(  S\xi\right)  _{n}=\sum_{k}a_{n-2k}\xi_{k}\text{\emph{;}}
\label{FourEq25}%
\end{equation}
and the adjoint $S^{\ast}$ is
\begin{equation}
\left(  S^{\ast}\xi\right)  _{n}=\sum_{k}\overline{a}_{k-2n}\xi_{k}.
\label{FourEq26}%
\end{equation}
\noindent\textup{(}b\textup{)} The 2D-dimensional subspace $\mathcal{L}$
spanned by $\left\{  e_{k}\mid-2D+1\leq k\leq0\right\}  $ is invariant under
$S^{\ast}$, and the corresponding $2D$ by $2D$ matrix is
\begin{equation}
\left(
\begin{array}
[c]{cccccccc}%
\bar{a}_{0} & 0 & 0 & \cdots & 0 & 0 & 0 & 0\\
\bar{a}_{2} & \bar{a}_{1} & \bar{a}_{0} & \cdots & \vdots & \vdots & \vdots &
\vdots\\
\vdots & \bar{a}_{3} & \bar{a}_{2} & \cdots & \vdots &  &  & \\
\bar{a}_{2D-2} & \vdots & \vdots &  & 0 & \vdots & \vdots & \vdots\\
0 & \bar{a}_{2D-1} & \bar{a}_{2D-2} & \ddots & \bar{a}_{0} & 0 & 0 & 0\\
0 & 0 & 0 &  & \bar{a}_{2} & \bar{a}_{1} & \bar{a}_{0} & 0\\
\vdots & \vdots & \vdots &  & \vdots & \bar{a}_{3} & \bar{a}_{2} & \bar{a}%
_{1}\\
\vdots & \vdots & \vdots & \cdots & \bar{a}_{2D-2} & \vdots & \vdots &
\vdots\\
0 & 0 & 0 & \cdots & 0 & \bar{a}_{2D-1} & \bar{a}_{2D-2} & \bar{a}_{2D-3}\\
0 & 0 & 0 & \cdots & 0 & 0 & 0 & \bar{a}_{2D-1}%
\end{array}
\right)  \text{.} \label{FourEq27}%
\end{equation}
Specifically, the case $D=2$ is
\begin{equation}
\left(
\begin{array}
[c]{cccc}%
\bar{a}_{0} & 0 & 0 & 0\\
\bar{a}_{2} & \bar{a}_{1} & \bar{a}_{0} & 0\\
0 & \bar{a}_{3} & \bar{a}_{2} & \bar{a}_{1}\\
0 & 0 & 0 & \bar{a}_{3}%
\end{array}
\right)  ; \label{FourEq28}%
\end{equation}
and the space $\mathcal{L}$ is then spanned by $\left\{  1,z^{-1}%
,z^{-2},z^{-3}\right\}  $.

\noindent\textup{(}c\textup{)} For all $n\in\mathbb{Z}$, there is a
$k\in\mathbb{N}$ such that
\begin{equation}
S^{\ast k}e_{n}\in\mathcal{L}\text{.} \label{FourEq29}%
\end{equation}

\end{lem}

\begin{pf}
The matrix representations \textup{(}\ref{FourEq25}\textup{)}, and
\textup{(}\ref{FourEq26}\textup{)}, follow directly from the formulas
\textup{(}\ref{FourEq23}\textup{)} and \textup{(}\ref{FourEq16}\textup{)}. To
prove \textup{(}b\textup{)}, set
\[
J\left(  D\right)  \text{:}=\left\{  k\in\mathbb{Z}\mid-2D+1\leq
k\leq0\right\}  \text{,}%
\]
so that
\[
\mathcal{L}=\operatorname*{span}\left\{  e_{k}\mid k\in J\left(  D\right)
\right\}  \text{.}%
\]
We claim the following implication:
\begin{equation}
\xi_{k}=0\text{ for }k\in\mathbb{Z}\backslash J\left(  D\right)
\Longrightarrow\sum_{k}\bar{a}_{k-2n}\xi_{k}=0\text{ for }n\in\mathbb{Z}%
\backslash J\left(  D\right)  \text{.} \label{FourEq30}%
\end{equation}
The invariance of $\mathcal{L}$ under $S^{\ast}$ is immediate from this. Now
suppose $\xi_{k}=0$ for $k\in\mathbb{Z}\backslash J\left(  D\right)  $. Then
one of the factors in $\bar{a}_{k-2n}\xi_{k}$ vanishes if $n\in\mathbb{Z}%
\backslash J\left(  D\right)  $; and the implication \textup{(}\ref{FourEq30}%
\textup{)} follows. To understand \textup{(}c\textup{)}, iterate the matrix
formula \textup{(}\ref{FourEq26}\textup{)}. We get
\begin{equation}
S^{\ast^{k}}e_{n}=\sum_{p}\left(  \sum_{i_{1},i_{2},\cdots,i_{k}}\bar
{a}_{i_{1}}\bar{a}_{i_{2}}\cdots\bar{a}_{i_{k}}\right)  e_{p} \label{FourEq31}%
\end{equation}
where the summation on the right-hand side is over $p\in\mathbb{Z}$,
$i_{1},i_{2},\cdots\in\lbrack0,1,\cdots,2D-1]$, subject to
\[
i_{1}+2i_{2}+\cdots+2^{k-1}i_{k}+p2^{k}=n\text{,}%
\]
so the range for the $p$-index is roughly divided by 2 with each iteration of
$S^{\ast}$ on the basis vector $e_{n}\left(  z\right)  =z^{n}$.
\end{pf}

\begin{rem}
\label{FourRem2}Conclusions \textup{(}a\textup{)}--\textup{(}b\textup{)} in
Lemma \ref{FourLem2} hold for the smaller subspace $\mathcal{M}$%
:$=\operatorname*{span}\left\{  e_{k}\mid-2D-2\leq k\leq0\right\}  $. The
matrix of the restricted operator $S^{\ast}\left\vert _{\overset{}%
{\mathcal{M}}}\right.  $ is then obtained from \textup{(}\ref{FourEq27}%
\textup{)} by deletion of the last row and last column. For example, the
reduced matrix corresponding to \textup{(}\ref{FourEq28}\textup{)} is
\[
F\text{\emph{:}}=\left(
\begin{array}
[c]{ccc}%
\bar{a}_{0} & 0 & 0\\
\bar{a}_{2} & \bar{a}_{1} & \bar{a}_{0}\\
0 & \bar{a}_{3} & \bar{a}_{2}%
\end{array}
\right)  \text{.}%
\]
However, property \textup{(}c\textup{)} in Lemma \ref{FourLem2} is not
satisfied for $\mathcal{M}$. To see this, note that $D=2$, and $\mathcal{L}%
=\mathcal{M}\oplus\mathbb{C}e_{-3}$ in the example. As a result
\begin{equation}
\left\langle e_{-3}\mid S^{\ast^{k}}e_{-3}\right\rangle =\left(  \bar{a}%
_{3}\right)  ^{k}\text{ for all }k\in\mathbb{N}\text{.} \label{FourEq33}%
\end{equation}
If $a_{3}\neq0$, it follows that $S^{\ast^{k}}e_{-3}$ is not in $\mathcal{M}$,
no matter how large $k$ is.
\end{rem}

\section{Measures Induced by Representations of the Cuntz
Relations\label{Measures}}

The system of coefficients studied in section \ref{FourPoly} give rise to
representations of the Cuntz relations. We recall the construction briefly;
see \cite{BrJo02b} for more details. When the coefficients $a_{0},a_{1}%
,\cdots$ are given, satisfying (\ref{FourEq3}), then the two operators $S_{0}$
and $S_{1}$ from Lemma \ref{FourLem1} (b) satisfy the relations (\ref{FourEq7}%
). We say that the two operators satisfy the Cuntz relations, or equivalently
that they define a representation of the Cuntz algebra $\mathcal{O}_{2}$ on
the Hilbert space $\mathcal{H}=L^{2}\left(  \mathbb{T}\right)  $. The
norm-closed algebra of operators on $\mathcal{H}$ generated by $S_{0}%
,S_{1},S_{0}^{\ast},$ and $S_{1}^{\ast}$ is known, \cite{Cu77} to be a simple
$C^{\ast}$-algebra; and it is unique up to isomorphism of $C^{\ast}$-algebras.

We will also need to Hilbert space $L^{2}\left(  \mathbb{R}\right)  $. If in
addition to (\ref{FourEq3}), the numbers $\left(  a_{k}\right)  $ satisfy
\begin{equation}
\sum_{k=0}^{2D-1}a_{k}=\sqrt{2}\text{,} \label{MeasEq1}%
\end{equation}
then it is known \cite{CMW92} that the following system of three equations has
solutions in $L^{2}\left(  \mathbb{R}\right)  $:
\begin{equation}
\varphi_{2n}\left(  x\right)  =\sqrt{2}\sum_{k}a_{k}\varphi_{n}\left(
2x-k\right)  \text{,} \label{MeasEq2}%
\end{equation}%
\begin{equation}
\varphi_{2n+1}\left(  x\right)  =\sqrt{2}\sum_{k}\left(  -1\right)
^{k}\overline{a_{2D-1-k}}\varphi_{n}\left(  2x-k\right)  \text{,}
\label{MeasEq3}%
\end{equation}
and
\begin{equation}
\int_{\mathbb{R}}\varphi_{0}\left(  x\right)  \;dx=1\text{.} \label{MeasEq4}%
\end{equation}
Moreover, for every $f\in L^{2}\left(  \mathbb{R}\right)  $, the following
\textit{Parseval identity} holds:
\begin{equation}
\sum_{n=0}^{\infty}\sum_{k\in\mathbb{Z}}\left\vert \left\langle \varphi
_{n}\left(  \cdot-k\right)  \mid f\right\rangle \right\vert ^{2}=\left\Vert
f\right\Vert ^{2}=\int_{\mathbb{R}}\left\vert f\left(  x\right)  \right\vert
^{2}\;dx\text{.} \label{MeasEq5}%
\end{equation}
When a system, such as
\begin{equation}
\left\{  \varphi_{n}\left(  \cdot-k\right)  \mid n\in\mathbb{N}_{0}%
,k\in\mathbb{Z}\right\}  \text{,} \label{MeasEq6}%
\end{equation}
satisfies (\ref{MeasEq5}) we way it is a \textit{Parseval frame}, or a
\textit{normalized tight frame}. There are simple further conditions on the
coefficients $\left(  a_{k}\right)  $ each of which guarantees that
(\ref{MeasEq6}) is in fact an orthonormal basis (ONB) in $L^{2}\left(
\mathbb{R}\right)  $, i.e., that
\begin{equation}
\left\langle \varphi_{n}\left(  \cdot-k\right)  \mid\varphi_{m}\left(
\cdot-\ell\right)  \right\rangle =\delta_{n,m}\delta_{k,\ell}\text{.}
\label{MeasEq7}%
\end{equation}

To understand how to select subsets $J\subset\mathbb{N}_{0}\times\mathbb{Z}$
such that the set of functions
\begin{equation}
\mathcal{F}\left(  J\right)  \text{:}=\left\{  2^{\frac{j}{2}}\varphi
_{n}\left(  2^{j}x-k\right)  \mid\left(  n,j\right)  \in J\text{, }%
k\in\mathbb{Z}\right\}  \label{MeasEq8}%
\end{equation}
forms a Perseval frame, or an ONB, for $L^{2}\left(  \mathbb{R}\right)  $, the
authors of \cite{CMW92} and \cite{Wic94} suggested a family of measures
associated with the operator systems (\ref{FourEq7}). We proved more
generally, in \cite{Jor04} that these measures are associated with any
representation of one of the Cuntz algebras $\mathcal{O}_{N}$, $N=2,3,\cdots$.

The idea is simple: Let $\left(  S_{i}\right)  _{i=0}^{1}$ be two operators in
a Hilbert space $\mathcal{H}$ such that
\begin{equation}
S_{i}^{\ast}S_{j}=\delta_{i,j}I\text{ and }\sum_{i=0}^{1}S_{i}S_{i}^{\ast}=I
\label{MeasEq9}%
\end{equation}
hold. Then for every multi-index
\begin{equation}
\xi\text{:}\left(  i_{1},i_{2},\cdots,i_{k}\right)  \text{, }i_{j}\in\left\{
0,1\right\}  \label{MeasEq10}%
\end{equation}
the operator
\begin{equation}
P\left(  \xi\right)  \text{:}=S_{i_{1}}\cdots S_{i_{k}}S_{i_{k}}^{\ast}\cdots
S_{i_{1}}^{\ast} \label{MeasEq11}%
\end{equation}
satisfies
\[
P\left(  \xi\right)  ^{\ast}=P\left(  \xi\right)  =P\left(  \xi\right)
^{2}\text{.}%
\]

Moreover
\begin{equation}
\sum_{\xi}P\left(  \xi\right)  =I \label{MeasEq12}%
\end{equation}
and
\begin{equation}
P\left(  \xi\right)  P\left(  \eta\right)  =\delta_{\xi,\eta}P\left(
\xi\right)  \label{MeasEq13}%
\end{equation}
where the summation in (\ref{MeasEq12}) is over all multi-indices of length
$k$, and in (\ref{MeasEq13}) the multi-indices $\xi$ and $\eta$ have the same
length. We say that, for each $k\in\mathbb{N}$ we have a dyadic partition of
the Hilbert space $\mathcal{H}$ into orthogonal \textquotedblleft frequency
bands.\textquotedblright

\begin{thm}
\label{MeasTheo0b}Let $m_{0}$\emph{:}$\mathbb{T\longrightarrow C}$ be a
function which satisfies $m_{0}\left(  1\right)  =\sqrt{2}$ and
\begin{equation}
\left\vert m_{0}\left(  z\right)  \right\vert ^{2}+\left\vert m_{0}\left(
-z\right)  \right\vert ^{2}=2\quad\text{for }z\in\mathbb{T}\text{.}
\label{MeasEq13a}%
\end{equation}
Set
\begin{equation}
m_{1}\left(  z\right)  =z\;\overline{m_{0}\left(  -z\right)  }\text{, for
}z\in\mathbb{T}\text{\emph{;}} \label{MeasEq13b}%
\end{equation}
and let $S_{i}$, $i=0,1$, be the operators \textup{(}\ref{FourEq6}\textup{)}.
Suppose the sequence $\varphi_{0},\varphi_{1},\varphi_{2},\cdots$ defined by
\begin{equation}
\widehat{\varphi_{2n}}\left(  x\right)  =\frac{1}{\sqrt{2}}m_{0}\left(
\frac{x}{2}\right)  \widehat{\varphi_{n}}\left(  \frac{x}{2}\right)
\label{MeasEq13c}%
\end{equation}%
\begin{equation}
\widehat{\varphi_{2n+1}}(x)=\frac{1}{\sqrt{2}}m_{1}\left(  \frac{x}{2}\right)
\widehat{\varphi_{n}}\left(  \frac{x}{2}\right)  \label{MeasEq13c1}%
\end{equation}
for $n=0,1,\cdots,\;x\in\mathbb{R}$ defines an orthonormal basis in
$L^{2}\left(  \mathbb{R}\right)  $, where we set $m_{j}\left(  x\right)
$\emph{:}$=m_{j}\left(  e^{i2\pi x}\right)  $, $j=0,1$. Then
\begin{equation}
2^{\frac{p}{2}}\varphi_{n}\left(  2^{p}t-k\right)  =\sum_{i_{1}=0}^{1}%
\cdots\sum_{i_{p}=0}^{1}\sum_{j\in\mathbb{Z}}\left\langle e_{j}\mid S_{i_{p}%
}^{\ast}\cdots S_{i_{1}}^{\ast}e_{k}\right\rangle \varphi_{2^{p}n+i_{1}%
+i_{2}2+\cdots+i_{p}2^{p-1}}\left(  t-j\right)  \label{MeasEq13d}%
\end{equation}
holds for all $p,n\in\mathbb{N}_{0}$ and $k\in\mathbb{Z}$. Moreover a subset
$J\subset\mathbb{N}_{0}\times\mathbb{N}_{0}$ has the corresponding set
\begin{equation}
\left\{  2^{\frac{p}{2}}\varphi_{n}\left(  2^{p}t-k\right)  \mid\left(
p,n\right)  \in J\text{, }k\in\mathbb{Z}\right\}  \label{MeasEq13e}%
\end{equation}
define an orthonormal basis \emph{(}ONB\emph{)} in $L^{2}\left(
\mathbb{R}\right)  $ if and only if the sets $\left[  2^{p}n,2^{p}\left(
n+1\right)  \right)  $ form a non-overlapping partition of $\mathbb{N}_{0}$.
\end{thm}

\begin{pf}
Since the functions in \textup{(}\ref{MeasEq13c}\textup{)}--\textup{(}%
\ref{MeasEq13c1}\textup{)} form an orthonormal basis for $L^{2}\left(
\mathbb{R}\right)  $, we have
\begin{equation}
\sum_{j\in\mathbb{Z}}\overline{\widehat{\varphi_{n}}\left(  x+j\right)
}\widehat{\varphi_{n^{\prime}}}\left(  x+j\right)  =\delta_{n,n^{\prime}%
}\text{ for }n,n^{\prime}\in\mathbb{N}_{0}\text{,} \label{MeasEq13f}%
\end{equation}
and \textup{(}\ref{MeasEq13d}\textup{)} will follow if we check that
\begin{equation}
\int_{\mathbb{R}}\overline{\varphi_{m}\left(  t-j\right)  }2^{\frac{p}{2}%
}\varphi_{n}\left(  2^{p}t-k\right)  \;dt=\left\langle e_{j}\mid S_{i_{p}%
}^{\ast}\cdots S_{i_{1}}^{\ast}e_{k}\right\rangle \label{MeasEq13g}%
\end{equation}
if $m=2^{p}n+i_{1}+i_{2}2+\cdots+i_{p}2^{p-1}$, and zero otherwise. But
\textup{(}\ref{MeasEq13g}\textup{)} follows from Fourier duality, and a
substitution of formulas \textup{(}\ref{MeasEq13c}\textup{)}--\textup{(}%
\ref{MeasEq13c1}\textup{)} and \textup{(}\ref{MeasEq13f}\textup{)}. We leave
the computation to the reader. The second conclusion regarding the basis
properties (ONB) of the functions in \textup{(}\ref{MeasEq13e}\textup{)}
follows from \textup{(}\ref{MeasEq13d}\textup{)} when we note that
\[
m=2^{p}n+i_{1}+i_{2}2+\cdots+i_{p}2^{p-1}%
\]
ranges over $\left[  2^{p}n,2^{p}\left(  n+1\right)  \right)  $ as
$i_{1},\cdots,i_{p}$ vary over $\underset{p\text{ times}}{\underbrace
{\mathbb{Z}_{2}\times\cdots\times\mathbb{Z}_{2}}}$.
\end{pf}

\begin{exmp}
\label{MeasEx0} The best known special case of \textup{(}\ref{MeasEq13e}%
\textup{)} is when
\[
J=\left\{  \left(  0,1\right)  ,\left(  1,1\right)  ,\left(  2,1\right)
,\left(  3,1\right)  ,\cdots\right\}  \text{,}%
\]
and the corresponding partition of $\mathbb{N}_{0}$ is simply $\left[
2^{p},2^{p+1}\right)  $ for $p=0,1,2,\cdots$. In this case, we only need the
first two functions from the sequence $\varphi_{n}$ of \textup{(}%
\ref{MeasEq13c}\textup{)}--\textup{(}\ref{MeasEq13c1}\textup{)}, i.e.,
$\varphi_{0}=\varphi$, and $\varphi_{1}=\psi$. The other extreme is
$J=\{(0,1),(0,2),(0,3),$ $\cdots\}$. In this case, each section in the
corresponding partition of $\mathbb{N}_{0}$ is a singleton, and \emph{all} the
functions $\left\{  \varphi_{n}\mid n\in\mathbb{N}_{0}\right\}  $ from
\textup{(}\ref{MeasEq13c}\textup{)}--\textup{(}\ref{MeasEq13c1}\textup{)} are
needed in the \textup{(}\ref{MeasEq13e}\textup{)} family. However, no dyadic
scales are used. Between these two extremes there is a countable infinite
family of non-overlapping partitions corresponding to the choice specified in
\textup{(}\ref{MeasEq13e}\textup{)}\emph{;} for example, the pairs $\left(
p,n\right)  $ in some $J$ from \textup{(}\ref{MeasEq13e}\textup{)} may be
specified as in the following table \textup{(}Table 1\textup{)}.
\[
\overline{\underline{%
\begin{array}
[c]{c}%
\text{\textbf{Table 1.}}\\
\text{A set }J\text{ consisting of points }\left(  p,n\right)  \text{ as
follows\emph{:}}\\%
\begin{array}
[c]{ccccccccccccccc}%
p & 2 & 2 & 2 & 2 & 4 & 4 & 4 & 6 & 6 & 6 & 8 & 8 & 8 & \cdots\\
n & 0 & 1 & 2 & 3 & 1 & 2 & 3 & 1 & 2 & 3 & 1 & 2 & 3 & \cdots
\end{array}
\end{array}
}}%
\]
In this case, we are using only the four functions $\varphi_{0},\varphi
_{1},\varphi_{2},$ and $\varphi_{3}$ from \textup{(}\ref{MeasEq13c}%
\textup{)}--\textup{(\ref{MeasEq13c1})}, and the segments from \textup{(}%
\ref{MeasEq13e}\textup{)} are $\left\{  0,1,2,3\right\}  $ followed by
\[
\left[  2^{2k}j,2^{2k}\left(  j+1\right)  \right)
\]
where $k=1,2,3,\cdots$, and $j\in\left\{  1,2,3\right\}  $. As a result, all
the dyadic scaling numbers $2^{2k}$ are used. The resulting ONB in
\textup{(}\ref{MeasEq13e}\textup{)} consists of the following functions,
$2\varphi_{0}\left(  2^{2}t-\ell_{0}\right)  ,$ $2^{k}\varphi_{j}\left(
2^{2k}t-\ell\right)  $, where $\ell_{0},\ell\in\mathbb{Z}$, $k=1,2,\cdots$,
and $j\in\left\{  1,2,3\right\}  $.
\end{exmp}

Since $\left\{  \varphi_{m}\left(  \cdot-j\right)  \mid m\in\mathbb{N}%
_{0},\text{ }j\in\mathbb{Z}\right\}  $ is an ONB for $L^{2}\left(
\mathbb{R}\right)  $, every $f\in L^{2}\left(  \mathbb{R}\right)  $,
$\left\Vert f\right\Vert ^{2}=1$, defines a probability distribution $P_{f}$
on $\mathbb{N}_{0}\times\mathbb{Z}$ by
\begin{equation}
P_{f}\left(  m,j\right)  \text{:}=\left\vert \left\langle \varphi_{m}\left(
\cdot-j\right)  \mid f\right\rangle \right\vert ^{2}\text{.} \label{MeasEq13h}%
\end{equation}

\begin{cor}
\label{MeasCor0} Let $p,n\in\mathbb{N}_{0}$, $k\in\mathbb{Z}$, and let
$P_{2^{\frac{p}{2}}\varphi_{n}\left(  2^{p}\cdot-k\right)  }\left(
\cdot,\cdot\right)  $ be the corresponding probability distribution from
\textup{(}\ref{MeasEq13h}\textup{)} above. Then the marginal distribution in
the first variable is
\begin{equation}
\delta_{m,2^{p}n+i_{1}+i_{2}2+\cdots+i_{p}2^{p-1}}\;\mu_{e_{k}}\left(  \left[
\frac{i_{1}}{2}+\cdots+\frac{i_{p}}{2^{p}},\frac{i_{1}}{2}+\cdots+\frac{i_{p}%
}{2^{p}}+\frac{1}{2^{p}}\right)  \right)  \text{.} \label{MeasEq13i}%
\end{equation}

\end{cor}

\begin{pf}
In view of \textup{(}\ref{MeasEq13h}\textup{)} and Theorem \ref{MeasTheo0b},
the marginal distribution is
\begin{align*}
&  \sum_{j\in\mathbb{Z}}P_{2^{\frac{p}{2}}\varphi_{n}\left(  2^{p}%
\cdot-k\right)  }\left(  m,j\right) \\
&  =\delta_{m,2^{p}n+i_{1}+i_{2}2+\cdots+i_{p}2^{p-1}}\sum_{j\in\mathbb{Z}%
}\left\vert \left\langle e_{j}\mid S_{i_{p}}^{\ast}\cdots S_{i_{1}}^{\ast
}e_{k}\right\rangle \right\vert ^{2}\\
&  =\delta_{m,2^{p}n+i_{1}+i_{2}2+\cdots+i_{p}2^{p-1}}\left\Vert S_{i_{p}%
}^{\ast}\cdots S_{i_{1}}^{\ast}e_{k}\right\Vert ^{2}\\
&  =\delta_{m,2^{p}n+i_{1}+i_{2}2+\cdots+i_{p}2^{p-1}}\;\mu_{e_{k}}\left(
\left[  \frac{i_{1}}{2}+\cdots+\frac{i_{p}}{2^{p}},\frac{i_{1}}{2}%
+\cdots+\frac{i_{p}}{2^{p}}+\frac{1}{2^{p}}\right)  \right)
\end{align*}
which is the desired conclusion.
\end{pf}

The next result follows from \cite{CMW92}, or from \cite{Jor04}.

\begin{lem}
\label{MeasLem1}There is a unique Borel measure $E$ defined on the unit
interval $\left[  0,1\right]  $ such that
\begin{equation}
E\left(  \left[  \frac{i_{1}}{2}+\cdots+\frac{i_{k}}{2^{k}},\frac{i_{1}}%
{2}+\cdots+\frac{i_{k}}{2^{k}}+\frac{1}{2^{k}}\right)  \right)  =P\left(
i_{1},i_{2},\cdots,i_{k}\right)  =S_{\xi}S_{\xi}^{\ast}\text{.}
\label{MeasEq14}%
\end{equation}
Moreover the measure $E=E_{P}$ takes values in orthogonal projections in
$\mathcal{H}$, and
\begin{equation}
\int_{0}^{1}E\left(  dx\right)  =I \label{MeasEq15}%
\end{equation}%
\begin{equation}
E\left(  B_{1}\right)  E\left(  B_{2}\right)  =0\text{ if }B_{i}\text{ are
Borel sets with }B_{1}\bigcap B_{2}=\varnothing\text{.} \label{MeasEq16}%
\end{equation}

\end{lem}

\begin{rem}
\label{MeasRem1}The expectation is that the measure class of $E$ is that of
the Lebesgue measure if $E$ is derived from a wavelet representation. For the
wavelet representations, the Hilbert space $\mathcal{H}$ is $L^{2}\left(
\mathbb{T}\right)  $; and for a generic set of wavelet representations $e_{0}$
is a cyclic vector; see \cite{Jor01}. Hence, to check if $E$ is in the
Lebesgue class, it is enough if we verify that the scalar valued, positive
measure
\begin{equation}
\mu_{0}\left(  \cdot\right)  \text{:}=\left\langle e_{0}\mid E\left(
\cdot\right)  e_{0}\right\rangle =\left\Vert E\left(  \cdot\right)
e_{0}\right\Vert ^{2} \label{MeasEq17}%
\end{equation}
is absolutely continuous with respect to Lebesgue measure on $\left[
0,1\right]  $. Using \textup{(}\ref{MeasEq14}\textup{)}, we \medskip see that
the measure $\mu_{0}$ is given on dyadic intervals by the formula
\begin{equation}
\mu_{0}\left(  \left[  \frac{i_{1}}{2}+\cdots+\frac{i_{k}}{2^{k}},\frac{i_{1}%
}{2}+\cdots+\frac{i_{k}}{2^{k}}+\frac{1}{2^{k}}\right)  \right)  =\left\Vert
S_{\xi}^{\ast}e_{0}\right\Vert ^{2}\text{.} \label{MeasEq18}%
\end{equation}

\end{rem}

Part of the reason for expecting that the measure $\mu_{0}\left(
\cdot\right)  =\left\Vert E\left(  \cdot\right)  e_{0}\right\Vert ^{2}$ should
be absolutely continuous with respect to Lebesgue measure on $\left[
0,1\right]  $ rests on the observation that $\lambda=\frac{1}{\sqrt{2}}$ is
always in the spectrum of the matrix \textup{(}\ref{FourEq27}\textup{)}. If
$\lambda=\frac{1}{\sqrt{2}}$ is also a dominant eigenvalue, in a sense which
we make precise below, then the asymptotic formula
\begin{equation}
\mu_{0}\left(  \left[  \frac{i_{1}}{2}+\cdots+\frac{i_{k}}{2^{k}},\frac{i_{1}%
}{2}+\cdots+\frac{1}{2^{k}}+\frac{1}{2^{k}}\right)  \right)  \approxeq2^{-k}
\label{MeasEq19}%
\end{equation}
would follow. From this the absolute continuity follows as well. But a closer
examination of the matrices $F_{i},i=0,1$, shows that $\lambda=\frac{1}%
{\sqrt{2}}$ is not always a dominant point in the spectrum of $F_{0}%
=S_{0}^{\ast}\left\vert _{\overset{}{\mathcal{M}}}\right.  $.

We now show that $\frac{1}{\sqrt{2}}$ is in the spectrum.

\begin{lem}
\label{MeasLem2}Let the numbers $a_{0},a_{1},\cdots,a_{2D-1}$ satisfy
\textup{(}\ref{FourEq3}\textup{)} and \textup{(}\ref{MeasEq1}\textup{)}, and
let $S_{0}$ and $S_{0}^{\ast}$ be the corresponding operators on $\ell^{2}$;
see \textup{(}\ref{FourEq25}\textup{)} and \textup{(}\ref{FourEq26}\textup{)}.
Then $\lambda=\frac{1}{\sqrt{2}}$ is an eigenvalue for the finite-dimensional
operator $F_{0}=S_{0}^{\ast}\left\vert _{\overset{}{\mathcal{M}}}%
\text{.}\right.  $ If $F_{0}$ is represented as a $2D-1$ by $2D-1$ matrix as
in \textup{(}\ref{FourEq27}\textup{)}, then the row-vector $w=\underset
{2D-1}{\underbrace{(1,1,\cdots,1)}}$ is a left-eigenvector for $F_{0}$, or
equivalently $w^{\ast}$ is a right-eigenvector, with the same eigenvalue
$\frac{1}{\sqrt{2}}$, for the adjoint matrix $F_{0}^{\ast}$.
\end{lem}

\begin{pf}
Substitution of $m_{0}\left(  1\right)  =\sqrt{2}$ into $\left\vert
m_{0}\left(  z\right)  \right\vert ^{2}+\left\vert m_{0}\left(  -z\right)
\right\vert ^{2}=2$ from Corollary \ref{FourCor1} yields $m_{0}\left(
-1\right)  =0$. Hence, $z+1$ divides the polynomial $m_{0}\left(  z\right)  $.
If
\begin{align*}
m_{0}\left(  z\right)   &  =\left(  z+1\right)  \mathcal{L}\left(  z\right)
\text{, and}\\
\mathcal{L}\left(  z\right)   &  =\sum_{i=0}^{2D-2}\ell_{i}z^{i}\text{, then
it follows that }\\
\mathcal{L}\left(  1\right)   &  =\sum_{i=0}^{2D-2}\ell_{i}=\frac{\sqrt{2}}%
{2}=\frac{1}{\sqrt{2}}\text{.}%
\end{align*}
As a result we get
\begin{equation}
\left\{
\begin{array}
[c]{l}%
a_{0}=\ell_{0}\\
a_{i}=\ell_{i}+\ell_{i-1}\text{ for }1\leq i\leq2D-2\\
a_{2D-1}=\ell_{2D-2}\text{.}%
\end{array}
\right.  \label{MeasEq20}%
\end{equation}
It follows from this that
\[
\sum_{i=0}^{D-1}a_{2i}=\sum_{i}\ell_{i}=\frac{1}{\sqrt{2}}%
\]
and
\[
\sum_{i=0}^{D-1}a_{2i+1}=\sum_{i}\ell_{i}=\frac{1}{\sqrt{2}}\text{.}%
\]
Working out the matrix product, this is a restatement of $wF_{0}=\frac
{1}{\sqrt{2}}w$, or equivalently $F_{0}^{\ast}w^{\ast}=\frac{1}{\sqrt{2}%
}w^{\ast}$.
\end{pf}

\begin{thm}
\label{MeasTheo0}Let $a_{0},a_{1},\cdots a_{2D-1}$ be as stated in Lemma
\ref{MeasLem2}, suppose
\[
\left\vert a_{0}\right\vert >\max\left\{  \left\vert \lambda\right\vert
\mid\lambda\in\operatorname*{spec}\left(  F_{0}\right)  \backslash\left\{
a_{0}\right\}  \right\}  \text{,}%
\]
and set
\[
s=-\frac{\ln\left\vert a_{0}\right\vert ^{2}}{\ln2}\text{.}%
\]
Then every non-empty open subset $V$ in $\left[  0,1\right]  $ contains an
infinite sequence of dyadic intervals $J$ such that
\[
0<\underset{J\subset V}{\lim\inf}\left(  \frac{\mu_{0}\left(  J\right)
}{\left\vert J\right\vert ^{s}}\right)  \leq\underset{J\subset V}{\lim\sup
}\left(  \frac{\mu_{0}\left(  J\right)  }{\left\vert J\right\vert ^{s}%
}\right)  <\infty\text{.}%
\]
We say that the measure $\mu_{0}$ contains $s$ as a fractal scale.
\end{thm}

\begin{pf}
The details follow by combining Lemma \ref{MeasLem2} with the technical lemma
which we state and prove in Section \ref{TechLem} below.
\end{pf}

\begin{prop}
\label{MeasProp1}Let the numbers $a_{0},\cdots,a_{2D-1}$ satisfy the two
conditions \textup{(}\ref{FourEq3}\textup{)} and \textup{(}\ref{MeasEq1}%
\textup{)} and let $S_{0},$ $S_{1}$ be the corresponding two operators; see
\textup{(}\ref{FourEq6}\textup{)}--\textup{(}\ref{FourEq7}\textup{)}. Let
$\mu_{0}\left(  \cdot\right)  $\emph{:}$=\left\Vert E\left(  \cdot\right)
e_{0}\right\Vert ^{2}$ be the measure \textup{(}\ref{MeasEq18}\textup{)}.
Then
\begin{equation}
\mu_{0}\left(  \left[  \frac{i_{1}}{2}+\cdots+\frac{i_{k}}{2^{k}},\frac{i_{1}%
}{2}+\cdots+\frac{i_{k}}{2^{k}}+\frac{1}{2^{k}}\right)  \right)
\geq\left\vert a_{0}\right\vert ^{2\cdot\#\left(  i=0\right)  }\left\vert
a_{2D-1}\right\vert ^{2\cdot\#\left(  i=1\right)  }\text{.} \label{MeasEq20a}%
\end{equation}

\end{prop}

\begin{rem}
\label{MeasRem2}We will see later how this estimate yields information about a
fractal component of $\mu_{0}$ in the special case $D=2$.
\end{rem}

\begin{pf}
Let $F_{i}$:$=S_{i}^{\ast}\left\vert \underset{\mathcal{M}}{}\right.  $,
$i=0,1$. For $k\in\mathbb{N}$, and $i_{1},\cdots,i_{k}\in\left\{  0,1\right\}
$, set $\xi=\frac{i_{1}}{2}+\cdots+\frac{i_{k}}{2^{k}}$ and $F_{\xi}%
$:$=F_{i_{k}}\cdots F_{i_{1}}$. Then $F_{0}^{\ast}e_{0}=\bar{a}_{0}e_{0}$
\noindent and $F_{1}^{\ast}e_{0}=\bar{a}_{2D-1}e_{0}$. To see this, recall
that
\[
m_{1}\left(  z\right)  =\bar{a}_{2D-1}-\bar{a}_{2D-2}z+\cdots-\bar{a}%
_{0}z^{2D-1}\text{,}%
\]
and
\begin{align*}
\mu_{0}\left(  \left[  \xi,\xi+2^{-k}\right)  \right)   &  =\left\Vert F_{\xi
}e_{0}\right\Vert ^{2}\geq\left\vert \left\langle e_{0}\mid F_{\xi}%
e_{0}\right\rangle \right\vert ^{2}\\
&  =\left\vert \left\langle F_{i_{1}}^{\ast}\cdots F_{i_{k}}^{\ast}e_{0}\mid
e_{0}\right\rangle \right\vert ^{2}\\
&  =\left\vert a_{0}\right\vert ^{2\cdot\#\left(  i=0\right)  }\left\vert
a_{2D-1}\right\vert ^{2\cdot\#\left(  i=1\right)  }\text{.}%
\end{align*}

\end{pf}

\begin{lem}
\label{MeasLem3}\textup{(}The case $D=2$.\textup{)}~

\noindent\textup{(}a\textup{)} For the case $D=2$, the real valued solutions
$a_{0},a_{1},a_{2},a_{3}$ to the three equations
\begin{equation}
\left\{
\begin{array}
[c]{l}%
a_{0}^{2}+a_{1}^{2}+a_{2}^{2}+a_{3}^{2}=1\\
a_{0}a_{2}+a_{1}a_{3}=0\\
a_{0}+a_{1}+a_{2}+a_{3}=\sqrt{2}%
\end{array}
\right.  \label{MeasEq21}%
\end{equation}
are
\begin{equation}
\left\{
\begin{array}
[c]{l}%
a_{0}=\frac{1}{2\sqrt{2}}\left(  1+\sqrt{2}\cos\beta\right) \\
a_{1}=\frac{1}{2\sqrt{2}}\left(  1+\sqrt{2}\sin\beta\right) \\
a_{2}=\frac{1}{2\sqrt{2}}\left(  1-\sqrt{2}\cos\beta\right) \\
a_{3}=\frac{1}{2\sqrt{2}}\left(  1-\sqrt{2}\sin\beta\right)
\end{array}
\right.  \label{MeasEq22}%
\end{equation}
for $\beta\in\mathbb{R}$.

\noindent\textup{(}b\textup{)}\emph{ }The spectrum of the matrix
\begin{equation}
F=\left(
\begin{array}
[c]{ccc}%
a_{0} & 0 & 0\\
a_{2} & a_{1} & a_{0}\\
0 & a_{3} & a_{2}%
\end{array}
\right)  \label{MeasEq23}%
\end{equation}
obtained by restriction of $S^{\ast}$ to $\mathcal{M}=\operatorname*{span}%
\left\{  e_{0},e_{-1},e_{-2}\right\}  $ is
\begin{equation}
\left\{  a_{0},\frac{1}{\sqrt{2}},\frac{\sin\beta-\cos\beta}{2}\right\}
\text{.} \label{MeasEq24}%
\end{equation}
\noindent\textup{(}c\textup{)} The following two conditions are equivalent
\begin{equation}
a_{0}>\frac{1}{\sqrt{2}}\geq\max\left\{  \left\vert \lambda\right\vert
\mid\lambda\in\operatorname*{spec}\left(  F\right)  \backslash\left\{
a_{0}\right\}  \right\}  \label{MeasEq25}%
\end{equation}
and
\begin{equation}
\left\vert \beta\right\vert <\frac{\pi}{4}\text{.} \label{MeasEq26}%
\end{equation}
\noindent\textup{(}d\textup{)} The eigenspace for the point $a_{0}$ in the
spectrum of $F$ is spanned by the vector $\binom{1}{v}$ where $v\in
\mathbb{R}^{2}$ is the solution to
\begin{equation}
v=\left(  a_{0}I_{2}-\left(
\begin{array}
[c]{cc}%
a_{1} & a_{0}\\
a_{3} & a_{2}%
\end{array}
\right)  \right)  ^{-1}\left(
\begin{array}
[c]{c}%
a_{2}\\
0
\end{array}
\right)  \text{.} \label{MeasEq27}%
\end{equation}
\noindent\textup{(}e\textup{)} For the case of Daubechies's wavelet, we have
\begin{equation}
\left\{
\begin{array}
[c]{l}%
a_{0}=\frac{1+\sqrt{3}}{4\sqrt{2}}\\
a_{1}=\frac{3+\sqrt{3}}{4\sqrt{2}}\\
a_{2}=\frac{3-\sqrt{3}}{4\sqrt{2}}\\
a_{3}=\frac{1-\sqrt{3}}{4\sqrt{2}}%
\end{array}
\text{,}\right.  \label{MeasEq28}%
\end{equation}%
\begin{equation}
\beta=\arccos\left(  \frac{\sqrt{3}-1}{2\sqrt{2}}\right)  \text{,}
\label{MeasEq29}%
\end{equation}
and
\begin{equation}
\operatorname*{spec}\left(  F\right)  =\left\{  a_{0},\frac{1}{\sqrt{2}}%
,\frac{1}{2\sqrt{2}}\right\}  \text{.} \label{MeasEq30}%
\end{equation}
\noindent In particular,
\begin{equation}
a_{0}<\frac{1}{\sqrt{2}}=\max\left\{  \left\vert \lambda\right\vert
\mid\lambda\in\operatorname*{spec}\left(  F\right)  \right\}  \text{.}
\label{MeasEq31}%
\end{equation}

\end{lem}

\begin{pf}
Part \textup{(}a\textup{)} follows directly from section 4 in
\cite{BraEvaJor2000}. To compute the spectrum of the matrix $F$, note that its
characteristic polynomial is
\begin{equation}
\left(  \lambda-a_{0}\right)  \left(  \lambda^{2}-\left(  a_{1}+a_{2}\right)
\lambda+a_{1}a_{2}-a_{0}a_{3}\right)  \text{;} \label{MeasEq32}%
\end{equation}
and when \textup{(}\ref{MeasEq22}\textup{)} are substituted, we see that the
roots are as listed in \textup{(}\ref{MeasEq24}\textup{)}. The conclusions
listed in parts \textup{(}c\textup{)} and \textup{(}d\textup{)} are immediate
consequences of \textup{(}b\textup{)}. Finally \textup{(}\ref{MeasEq27}%
\textup{)} follows from the observation that $F$ has the form
\begin{equation}
F=\left(
\begin{array}
[c]{cc}%
a_{0} &
\begin{array}
[c]{cc}%
0\; & 0
\end{array}
\\%
\begin{array}
[c]{c}%
a_{2}\\
0
\end{array}
& G
\end{array}
\right)  \label{MeasEq33}%
\end{equation}
with $G=\left(
\begin{array}
[c]{cc}%
a_{1} & a_{0}\\
a_{3} & a_{2}%
\end{array}
\right)  $ and the fact that $a_{0}$ is not in the spectrum of $G$. The
conclusions listed in \textup{(}e\textup{)} for Daubechies's wavelet follow by
an application of \textup{(}a\textup{)}--\textup{(}d\textup{)}.
\end{pf}

To better appreciate the geometry of the formulas (\ref{MeasEq22}), the reader
may find Fig. 1 useful.

\[%
{\parbox[b]{4.5383in}{\begin{center}
\includegraphics[
natheight=6.985600in,
natwidth=14.145400in,
height=2.2506in,
width=4.5383in
]%
{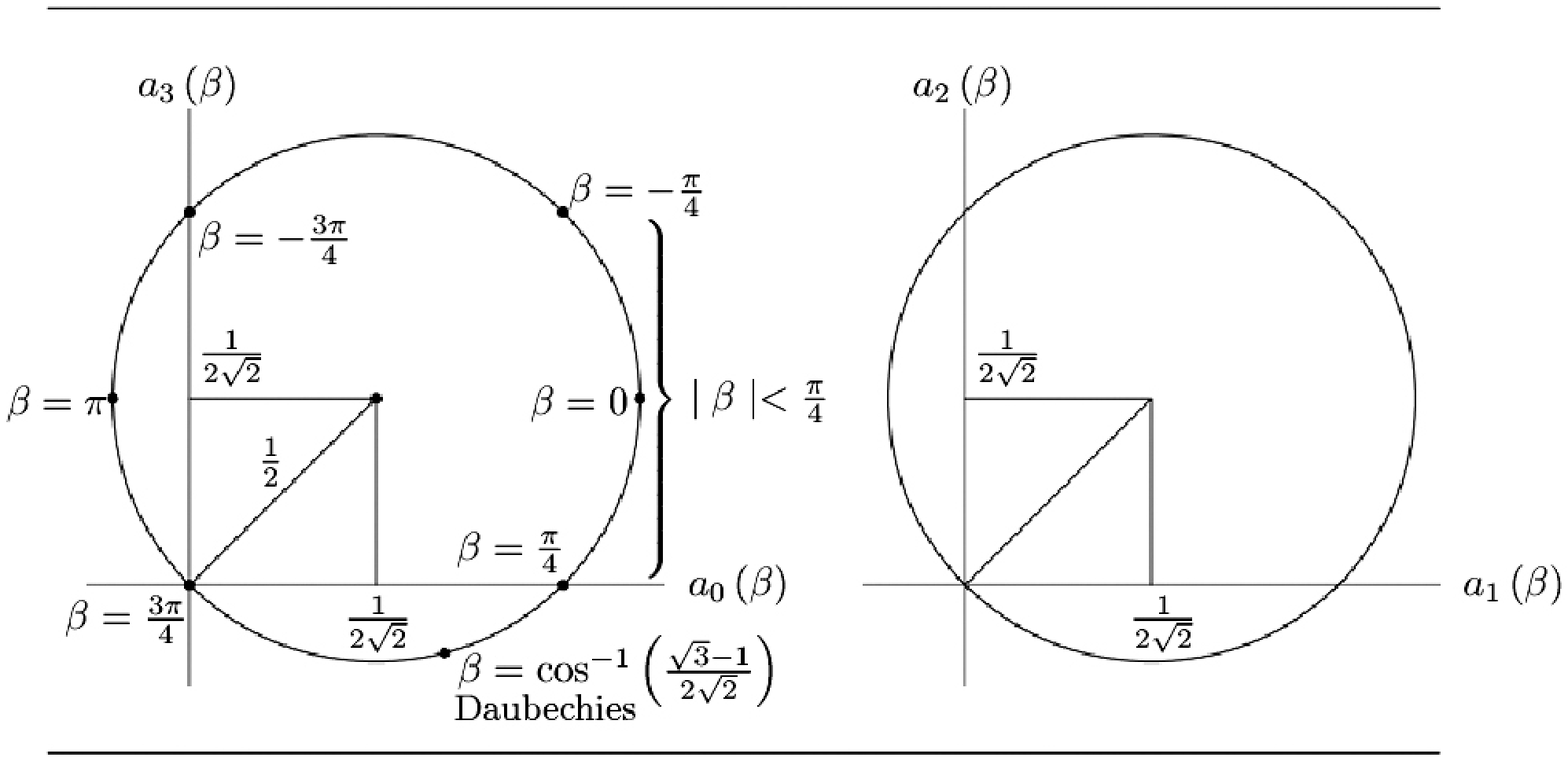}%
\\
The coefficients $a_0\left(  \beta\right)  ,\,a_1\left(  \beta\right)
,\,a_2\left(  \beta\right)  ,\,a_3\left(  \beta\right)  $ of (\ref{MeasEq21}%
)\medskip\newline Figure 1.
\end{center}}}
\]

The four coefficients $a_{0},a_{1},a_{2},a_{3}$ are described by the two
circles in Fig. 1, viz., by
\begin{equation}
\left(  a_{0}-\frac{1}{2\sqrt{2}}\right)  ^{2}+\left(  a_{3}-\frac{1}%
{2\sqrt{2}}\right)  ^{2}=\frac{1}{4}\text{, and by }\left(  a_{1}-\frac
{1}{2\sqrt{2}}\right)  ^{2}+\left(  a_{2}-\frac{1}{2\sqrt{2}}\right)
^{2}=\frac{1}{4}\text{.} \label{MeasEq34}%
\end{equation}
This representation also makes it clear that, if one of the four coefficients
vanishes, then so does a second from the remaining coefficients. The resulting
four degenerate cases are as sketched in Table 2.\bigskip\
\[
\overline{\underline{%
\begin{array}
[c]{c}%
\text{\textbf{\underline{Table\ 2.}}}\\
\multicolumn{1}{l}{\text{{\small The degenerate cases when one (and
therefore}}}\\
\multicolumn{1}{l}{\text{{\small two) of the coefficients} }a_{0},a_{1}%
,a_{2},a_{3}\text{ {\small vanish:}}}\\
\\%
\begin{array}
[c]{ccccc}%
m_{0} & \frac{1+z}{\sqrt{2}} & \frac{1+z^{3}}{\sqrt{2}} & \frac{z+z^{2}}%
{\sqrt{2}} & \frac{z^{2}+z^{3}}{\sqrt{2}}\\
&  &  &  & \\
m_{1} & \frac{z^{2}-z^{3}}{\sqrt{2}} & \frac{1-z^{3}}{\sqrt{2}} &
\frac{-z+z^{2}}{\sqrt{2}} & \frac{1-z}{\sqrt{2}}\\
&  &  &  & \\
\beta & \frac{\pi}{4} & -\frac{\pi}{4} & \frac{3\pi}{4} & -\frac{3\pi}{4}%
\end{array}
\end{array}
}}%
\]

\bigskip

\noindent There are four distinct configurations, and all four are variations
of the Haar wavelet, see Fig. 2. The second column in Table 2 is called the
stretched Haar wavelet; see \cite{BrJo02b}. For all three, we have the formula%

\begin{equation}
\mu_{0}\left(  \left[  \frac{i_{1}}{2}+\cdots+\frac{i_{k}}{2^{k}},\frac{i_{1}%
}{2}+\cdots+\frac{i_{k}}{2^{k}}+\frac{1}{2^{k}}\right)  \right)  =2^{-k}
\label{MeasEq35}%
\end{equation}
which proves that $\mu_{0}$ \textit{is} the Lebesgue measure on the
unit-interval $\left[  0,1\right]  $ for all four variations of the Haar
wavelet. The two matrices $F_{0}$ and $F_{1}$ corresponding to the four
variants of the Haar wavelet filters in Table 2 are as follows:

\vspace{0.5in}

\noindent$%
\begin{array}
[c]{l}%
F_{0}\text{:\quad}\frac{1}{\sqrt{2}}\left(
\begin{array}
[c]{ccc}%
1 & 0 & 0\\
0 & 1 & 1\\
0 & 0 & 0
\end{array}
\right)  \frac{1}{\sqrt{2}}\left(
\begin{array}
[c]{ccc}%
1 & 0 & 0\\
0 & 0 & 1\\
0 & 1 & 0
\end{array}
\right)  \frac{1}{\sqrt{2}}\left(
\begin{array}
[c]{ccc}%
0 & 0 & 0\\
1 & 1 & 0\\
0 & 0 & 1
\end{array}
\right)  \frac{1}{\sqrt{2}}\left(
\begin{array}
[c]{ccc}%
0 & 0 & 0\\
1 & 0 & 0\\
0 & 1 & 1
\end{array}
\right) \\
\\
F_{1}\text{:\quad}\frac{1}{\sqrt{2}}\!\left(
\begin{array}
[c]{ccc}%
\text{0} & \text{0} & \text{0}\\
\text{1} & \text{0} & \text{0}\\
\text{0} & \text{-1} & \text{1}%
\end{array}
\right)  \frac{1}{\sqrt{2}}\!\left(
\begin{array}
[c]{ccc}%
\text{1} & \text{0} & \text{0}\\
\text{0} & \text{0} & \text{1}\\
\text{0} & \text{-1} & \text{0}%
\end{array}
\right)  \frac{1}{\sqrt{2}}\!\left(
\begin{array}
[c]{ccc}%
\text{0} & \text{0} & \text{0}\\
\text{1} & \text{-1} & \text{0}\\
\text{0} & \text{0} & \text{1}%
\end{array}
\right)  \frac{1}{\sqrt{2}}\!\left(
\begin{array}
[c]{ccc}%
\text{1} & \text{0} & \text{0}\\
\text{0} & \text{-1} & \text{1}\\
\text{0} & \text{0} & \text{0}%
\end{array}
\right)
\end{array}
$ \newpage%

\[
\underline{\overline{%
\begin{array}
[c]{c}%
\text{The father function }\varphi\text{ and the mother function }\psi\text{
for the Haar Wavelet}\\
\text{%
{\parbox[b]{4.0491in}{\begin{center}
\includegraphics[
natheight=2.794800in,
natwidth=4.173700in,
height=2.719in,
width=4.0491in
]%
{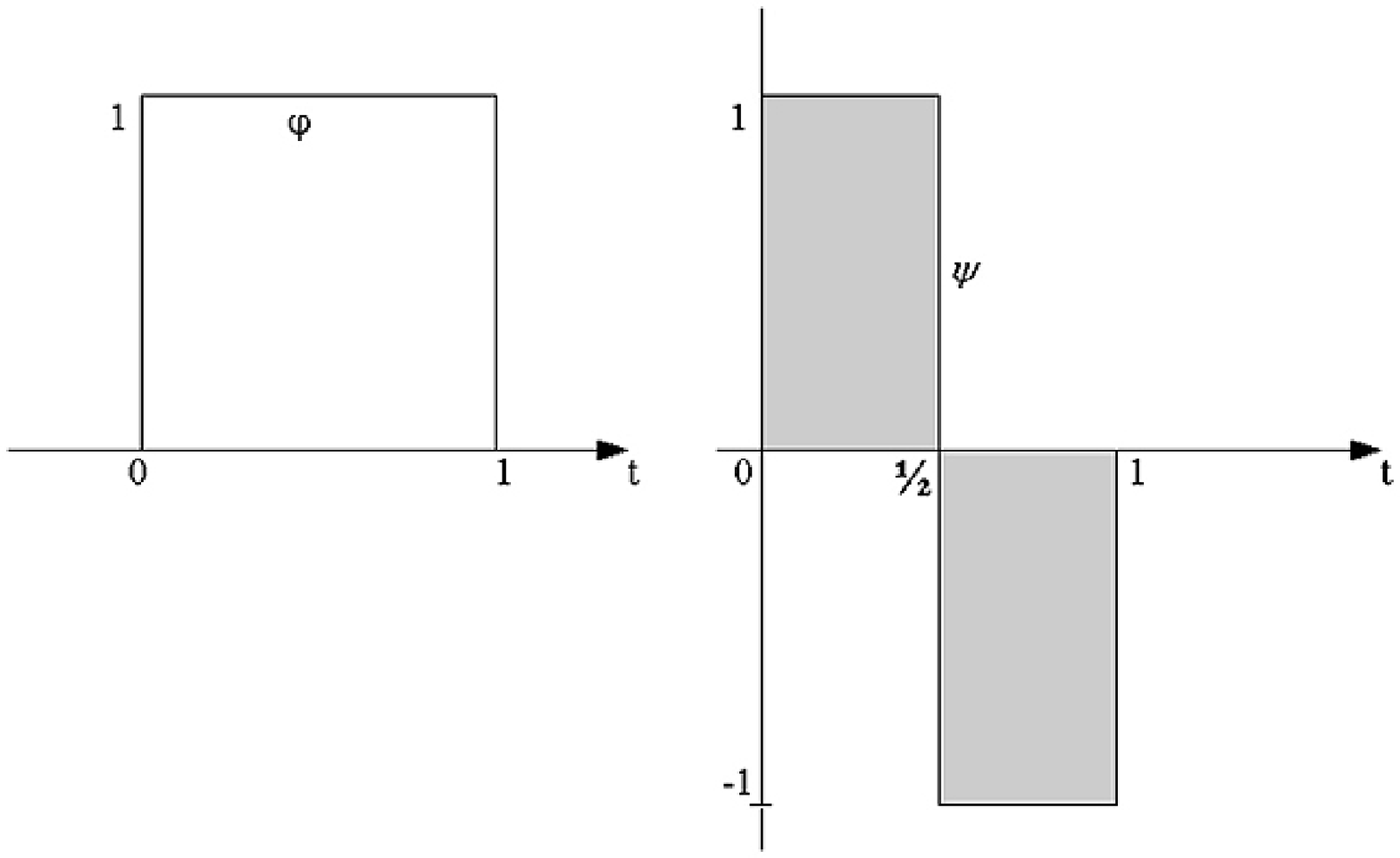}%
\\
Figure 2.
\end{center}}}
}%
\end{array}
}}%
\]

\begin{rem}
\label{MeasRem3}We shall need Lemma \ref{MeasLem2} for the two matrices
\begin{equation}
F_{0}=\left(
\begin{array}
[c]{ccc}%
a_{0} & 0 & 0\\
a_{2} & a_{1} & a_{0}\\
0 & a_{3} & a_{2}%
\end{array}
\right)  \text{ and }F_{1}=\left(
\begin{array}
[c]{ccc}%
a_{3} & 0 & 0\\
a_{1} & -a_{2} & a_{3}\\
0 & -a_{0} & a_{3}%
\end{array}
\right)  \label{MeasEq36}%
\end{equation}
For the range $\left\vert \beta\right\vert <\frac{\pi}{4}$, we have
\textup{(}\ref{MeasEq25}\textup{)}\emph{ }satisfied by $\operatorname*{spec}%
\left(  F_{0}\right)  $. But then the point $a_{3}=a_{3}\left(  \beta\right)
$ satisfies $\left\vert a_{3}\right\vert <\frac{1}{\sqrt{2}}$, and
\begin{equation}
\frac{1}{\sqrt{2}}=\max\left\{  \left\vert \lambda\right\vert \mid\lambda
\in\operatorname*{spec}\left(  F_{1}\right)  \right\}  \text{.}
\label{MeasEq37}%
\end{equation}

\end{rem}

\begin{thm}
\label{MeasTheo1}Let $\beta\in\mathbb{R}$, $\left\vert \beta\right\vert
<\frac{\pi}{4}$, be given. Let $a_{0},$ $a_{1},$ $a_{2},$ $a_{3},$ be the
corresponding numbers in \textup{(}\ref{MeasEq22}\textup{)}. Let $\left(
S_{i}\right)  _{i=0}^{1}$ be the operators in \textup{(}\ref{FourEq6}%
\textup{)}, and
\begin{equation}
F_{i}\text{:}=S_{i}^{\ast}\left\vert \underset{\mathcal{M}}{}\right.
\label{MeasEq38}%
\end{equation}
the matrices \textup{(}\ref{MeasEq36}\textup{)}. Let
\begin{equation}
\xi=\frac{i_{1}}{2}+\cdots+\frac{i_{k}}{2^{k}} \label{MeasEq39}%
\end{equation}
be a dyadic fraction. Then
\begin{equation}
\lim_{n\longrightarrow\infty}a_{0}^{-n}F_{0}^{n}F_{i_{k}}\cdots F_{i_{1}}%
e_{0}=a_{0}^{\#\left(  i=0\right)  }a_{3}^{\#\left(  i=1\right)  }\left(
e_{0}+v\right)  \label{MeasEq40}%
\end{equation}
where $v$ is the vector in \textup{(}\ref{MeasEq27}\textup{)}. Moreover, if
$\varepsilon>0$ is given, there is an $n_{0}\in\mathbb{N}$ such that
\begin{equation}
\mu_{0}\left(  \left[  \xi,\xi+\frac{1}{2^{n+k}}\right)  \right)  \geq
a_{0}^{2n}a_{0}^{2\cdot\#\left(  i=0\right)  }a_{3}^{2\cdot\#\left(
i=1\right)  }\left(  1+\left\Vert v\right\Vert ^{2}-\varepsilon\right)
\text{.} \label{MeasEq41}%
\end{equation}
for all $n\geq n_{0}$. More generally, given $\varepsilon$, $0<\varepsilon<1$,
there is $n_{0}\in\mathbb{N}$ such that
\begin{align*}
\left(  1+\left\Vert v\right\Vert ^{2}-\varepsilon\right)  \left(
a_{0}^{\#\left(  i=0\right)  }a_{3}^{\#\left(  i=1\right)  }\right)  ^{2}  &
\leq\frac{\mu_{0}\left(  \left[  \xi,\xi+\frac{1}{2^{n+k}}\right)  \right)
}{a_{0}^{2n}}\\
&  \leq\left(  a_{0}^{\#\left(  i=0\right)  }a_{3}^{\#\left(  i=1\right)
}\right)  ^{2}\left(  1+\left\Vert v\right\Vert ^{2}+\varepsilon\right)
\text{ for all }n\geq n_{0}\text{.}%
\end{align*}

\end{thm}

\begin{pf}
Recall that $e_{0}+v$ is the eigenvector satisfying
\begin{equation}
F_{0}\left(  e_{0}+v\right)  =a_{0}\left(  e_{0}+v\right)  \text{, and
}\left\langle e_{0}\mid e_{0}+v\right\rangle =1\text{.} \label{MeasEq42}%
\end{equation}
Moreover, we have
\begin{equation}
F_{0}^{\ast}e_{0}=a_{0}e_{0}\text{ and }F_{1}^{\ast}e_{0}=a_{3}e_{0}\text{.}
\label{MeasEq43}%
\end{equation}
Using \textup{(}\ref{MeasEq25}\textup{)} for $F_{0}$, we get
\begin{equation}
\lim_{n\longrightarrow\infty}a_{0}^{-n}F_{0}^{n}w=\left\langle e_{0}\mid
w\right\rangle \left(  e_{0}+v\right)  \label{MeasEq44}%
\end{equation}
for all $w\in\mathcal{M}$. Applying this, and \textup{(}\ref{MeasEq43}%
\textup{)}, to $w$:$=F_{i_{k}}\cdots F_{i_{1}}e_{0}$, the conclusion
\textup{(}\ref{MeasEq40}\textup{)} follows. Since
\begin{equation}
\mu_{0}\left(  \left[  \xi,\xi+\frac{1}{2^{n+k}}\right)  \right)  =\left\Vert
S_{0}^{\ast^{n}}S_{\xi}^{\ast}e_{0}\right\Vert ^{2}=\left\Vert F_{0}%
^{n}F_{i_{k}}\cdots F_{i_{1}}e_{0}\right\Vert ^{2}\text{,} \label{MeasEq45}%
\end{equation}
the second conclusion \textup{(}\ref{MeasEq41}\textup{)} follows from the first.
\end{pf}

\begin{cor}
\label{MeasCor1}Let $\left\vert \beta\right\vert <\frac{\pi}{4}$, and let
$a_{0},$ $a_{1},$ $a_{2},$ $a_{3}$ be the corresponding numbers in
\textup{(}\ref{MeasEq22}\textup{)}. If the measure $\mu_{0}\left(
\cdot\right)  =\left\Vert E\left(  \cdot\right)  e_{0}\right\Vert ^{2}$ is
absolutely continuous with respect to Lebesgue measure $dx$ on $\left[
0,1\right]  $. Then the Radon-Nikodym derivative $d\mu_{0}/dx=f$ is unbounded
in every open \textup{(}non-empty\textup{)} subset of $\left[  0,1\right]  $.
\end{cor}

\begin{pf}
Let $V\subset\left[  0,1\right]  $ be a non-empty open subset. Then pick $\xi$
as in \textup{(}\ref{MeasEq39}\textup{)}, and $n_{1}\in\mathbb{N}$ such that
\begin{equation}
\left[  \xi,\xi+\frac{1}{2^{n_{1}}}\right)  \subset V\text{.} \label{MeasEq46}%
\end{equation}
Then pick $\varepsilon$, $0<\varepsilon<\frac{1}{2}$, and $n_{0}\geq n_{1}$
such that \textup{(}\ref{MeasEq41}\textup{)} holds for all $n\geq n_{0}$. If
$f$ were bounded in $V$ with upper bound $C$; then \textup{(}\ref{MeasEq41}%
\textup{)} implies the estimate
\begin{equation}
C\cdot2^{-n}\geq a_{0}^{2n}a_{0}^{2\cdot\#\left(  i=0\right)  }a_{3}%
^{2\cdot\#\left(  i=1\right)  }\left(  1+\left\Vert v\right\Vert
^{2}-\varepsilon\right)  \label{MeasEq47}%
\end{equation}
for all $n\geq n_{0}$. Since $a_{0}^{2}>\frac{1}{2}$, $\lim_{n\longrightarrow
\infty}2^{n}a_{0}^{2n}=\infty$, which contradicts \textup{(}\ref{MeasEq47}%
\textup{)}.
\end{pf}

\begin{exmp}
\label{MeasEx1}We introduced the constants $a_{0},$ $a_{1},$ $a_{2},$ $a_{3}$
for Daubechies's wavelet in part \textup{(}e\textup{)} of Lemma \ref{MeasLem2}%
. An inspection of \textup{(}\ref{MeasEq22}\textup{)} and \textup{(}%
\ref{MeasEq28}\textup{)} shows that $\beta=\arccos\left(  \frac{\sqrt{3}%
-1}{2\sqrt{2}}\right)  $. Hence
\[
\frac{\sin\beta-\cos\beta}{2}=\frac{1}{2\sqrt{2}}\text{,}%
\]
and the three numbers in the spectrum of $F_{0}$ are as follows\emph{:}
\begin{equation}
\frac{1}{2\sqrt{2}}<a_{0}<\frac{1}{\sqrt{2}}\text{.} \label{MeasEq48}%
\end{equation}
One checks that the vector
\begin{equation}
v^{\ast}\text{\emph{:}}=\left(  1+\sqrt{3}\right)  e_{-1}+\left(  1-\sqrt
{3}\right)  e_{-2} \label{MeasEq49}%
\end{equation}
spans the $\frac{1}{\sqrt{2}}$-eigenspace for $F_{0}$. Since $a_{0}<\frac
{1}{\sqrt{2}}$, the eigenspace
\begin{equation}
\left\{  w\in\mathcal{M}\mid F_{0}w=\frac{1}{\sqrt{2}}w\right\}
\label{MeasEq50}%
\end{equation}
is orthogonal to $e_{0}$.
\end{exmp}

\begin{cor}
\label{MeasCor2}Let $\beta=\arccos\left(  \frac{\sqrt{3}-1}{2\sqrt{2}}\right)
$, so that the numbers $a_{0},$ $a_{1},$ $a_{2},$ $a_{3}$ define the
Daubechies wavelet. Let $v^{\ast}$ be the eigenvector in \textup{(}%
\ref{MeasEq49}\textup{)}, and pick $w\in\mathcal{M}$ such that
\begin{equation}
F_{0}^{\ast}w=\frac{1}{\sqrt{2}}w\text{, and }\left\langle w\mid v^{\ast
}\right\rangle =\sqrt{8}\text{.} \label{MeasEq51}%
\end{equation}
Let $\xi$ be a dyadic rational given as in \textup{(}\ref{MeasEq39}\textup{)}.
Then, for every $\varepsilon$, $0<\varepsilon<1$, there is $n_{0}\in
\mathbb{N}$ such that
\begin{align}
\left\vert \left\langle w\mid F_{i_{k}}\cdots F_{i_{1}}e_{0}\right\rangle
\right\vert ^{2}-\varepsilon &  \leq\frac{\mu_{0}\left(  \left[  \xi,\xi
+\frac{1}{2^{n+k}}\right)  \right)  }{2^{-n}}\label{MeasEq52}\\
&  \leq\left\vert \left\langle w\mid F_{i_{k}}\cdots F_{i_{1}}e_{0}%
\right\rangle \right\vert ^{2}+\varepsilon\text{ for all }n\geq n_{0}%
\text{.}\nonumber
\end{align}

\end{cor}

\begin{pf}
The argument follows the proof of Theorem \ref{MeasTheo1}. The vector $w$ is
$w=\sqrt{2}\left(  e_{0}+e_{-1}+e_{-2}\right)  $. The main difference is that
now $\frac{1}{\sqrt{2}}$ is the top eigenvalue of $F_{0}$ with eigenvector
$v^{\ast}$. If $w$ is chosen as in \textup{(}\ref{MeasEq51}\textup{)}, then
\begin{equation}
\lim_{n\longrightarrow\infty}2^{\frac{n}{2}}F_{0}^{n}s=\frac{1}{\sqrt{8}%
}\left\langle w\mid s\right\rangle v^{\ast} \label{MeasEq53}%
\end{equation}
for all $s\in\mathcal{M}$. The desired result \textup{(}\ref{MeasEq52}%
\textup{)} follows when this is applied to $s=F_{i_{k}}\cdots F_{i_{1}}e_{0}$.
\end{pf}

\section{Wavelets on Fractals\label{WaveFract}}

In section \ref{FourPoly} we proved some lemmas for representations of the
Cuntz algebra $\mathcal{O}_{2}$, and we used them in section \ref{Measures}
above in our analysis of dyadic wavelets in the Hilbert space $L^{2}\left(
\mathbb{R}\right)  $. But there are other Hilbert spaces $\mathcal{H}$ which
admit wavelet algorithms. In a recent paper \cite{DutJor04}, we constructed
wavelets in separable Hilbert spaces $\mathcal{H}^{(s)}$ built on Hausdorff
measure $(dx)^{s}$ where $s$ denotes the corresponding fractal dimension. For
the middle third Cantor set, for example, the scaling number $N$ is $3$; and
the fractal dimension is $s=\frac{\ln2}{\ln3}=\log_{3}(2)$.

In this section, we show how our Cantor subdivision construction is built on a
representation of the $\mathcal{O}_{3}$-Cuntz relations on $L^{2}\left(
\mathbb{T}\right)  $ which are analogous to the representations of
$\mathcal{O}_{2}$ which we used in our analysis of dyadic $L^{2}\left(
\mathbb{R}\right)  $-wavelets.

\begin{lem}
\label{WaveFractLem1}Let $N\in\mathbb{N}$, $N\geq2$, be given, and let
$m_{0},$ $m_{1},\cdots,m_{N-1}$ be bounded measurable functions on
$\mathbb{T}$. Then the operators
\begin{equation}
\left(  S_{i}f\right)  \left(  z\right)  =m_{i}\left(  z\right)  f\left(
z^{N}\right)  \text{, }f\in L^{2}\left(  \mathbb{T}\right)  \text{,
}i=0,1,\cdots,N-1\text{, }z\in\mathbb{T} \label{WavFracEq1}%
\end{equation}
satisfy the Cuntz relations
\begin{equation}
\left\{
\begin{array}
[c]{c}%
S_{i}^{\ast}S_{j}=S_{i,j}I\\
\sum\limits_{i=0}^{N-1}S_{i}S_{i}^{\ast}=I
\end{array}
\right.  \label{WavFracEq2}%
\end{equation}
if and only if the associated $N$ by $N$ matrix function
\begin{equation}
M_{N}\left(  z\right)  =\frac{1}{\sqrt{N}}\left(  m_{j}\left(  ze^{i\frac
{k2\pi}{N}}\right)  \right)  _{j,k} \label{WavFracEq3}%
\end{equation}
takes values in the unitary matrices for a.e. $z\in\mathbb{T}$.
\end{lem}

\begin{pf}
Rather than sketching the details, we will instead refer the reader to section
\ref{FourPoly} above where the argument is done in full for $N=2$; see also
\cite{DutJor04} or \cite{BrJo02b}.\qed

\end{pf}

For each of the representations $\left(  S_{i}\right)  _{i=0}^{N-1}$, we get a
measure $\mu_{0}$ on the unit-interval $\left[  0,1\right]  $, as outlined in
section \ref{Measures} (for the special case $N=2$.) In the general case, the
measure $\mu_{0}$ is determined uniquely on the $N$-adic subintervals as
follows,
\begin{equation}
\mu_{0}\left(  \left[  \frac{i_{1}}{N}+\cdots+\frac{i_{k}}{N^{k}}\text{,
}\frac{i_{1}}{N}+\cdots+\frac{i_{k}}{N^{k}}+\frac{1}{N^{k}}\right)  \right)
=\left\Vert S_{i_{k}}^{\ast}\cdots S_{i_{2}}^{\ast}S_{i_{1}}^{\ast}%
e_{0}\right\Vert ^{2} \label{WavFracEq4}%
\end{equation}
where $i_{1},i_{2},\cdots\in\left\{  0,1,\cdots,N-1\right\}  $.

We state the next result just for the middle-third Cantor set; but, following
the discussion in \cite{DutJor04}, the reader will convince him/herself that
it carries over to any fractal constructed by iteration of finite families of
(contractive) affine maps in $\mathbb{R}^{d}$; so called iterated function
systems (IFS).

\begin{prop}
\label{WaveFractProp1}Let
\begin{equation}
\left\{
\begin{array}
[c]{l}%
m_{0}(z)=\frac{1+z^{2}}{\sqrt{2}}\\
m_{1}\left(  z\right)  =z\\
m_{2}\left(  z\right)  =\frac{1-z^{2}}{\sqrt{2}}%
\end{array}
\right.  \label{WavFracEq5}%
\end{equation}
Then the unitarity condition \textup{(}\ref{WavFracEq3}\textup{)} from Lemma
\ref{WaveFractLem1} is satisfied, and the operators
\begin{equation}
\left(  S_{i}f\right)  \left(  z\right)  =m_{i}\left(  z\right)  f\left(
z^{3}\right)  \text{, }f\in L^{2}\left(  \mathbb{T}\right)  \text{,
}i=0,1,2\text{, }z\in\mathbb{T}\text{,} \label{WavFracEq6}%
\end{equation}
define a representation of $\mathcal{O}_{3}$ on $L^{2}\left(  \mathbb{T}%
\right)  $; and the corresponding measure $\mu_{0}$ from \textup{(}%
\ref{WavFracEq4}\textup{)} is the Hausdorff measure $\mu^{\left(  s\right)  }$
of Hausdorff dimension $s=\frac{\ln2}{\ln3}=\log_{3}\left(  2\right)  $
restricted to the middle-third Cantor set $X_{3}\subset\left[  0,1\right]  $.
\end{prop}

\begin{rem}
\label{WaveFractRem1}A classical theorem of Hutchinson \cite{Hut81}, see also
\cite{Falconer85}, implies that the Cantor measure, supported on $X_{3},$ is
the unique Borel probability measure $\mu$ which solves
\begin{align}
\int f\left(  x\right)  \;dx\left(  x\right)   &  =\frac{1}{2}\left(  \int
f\left(  \frac{x}{3}\right)  \;dx\left(  x\right)  +\int f\left(  \frac
{x+2}{3}\right)  \;d\mu\left(  x\right)  \right) \label{WavFracEq7}\\
&  \text{for all continuous bounded functions }f\text{.}\nonumber
\end{align}
As a result, we need to verify that our measure $\mu_{0}$ from \textup{(}%
\ref{WavFracEq4}\textup{)}, and from the $\mathcal{O}_{3}$-representation,
satisfies identity \textup{(}\ref{WavFracEq7}\textup{)}.
\end{rem}

To compute the terms on the right-hand side in \textup{(}\ref{WavFracEq4}%
\textup{)} is a finite matrix problem. We need only to calculate the operators
$S_{i}^{\ast}$, $i=0,1,2,$ on the three-dimensional subspace $\mathcal{M}%
$:$=\operatorname*{span}\left\{  e_{0},e_{-1},e_{-2}\right\}  $. But the
argument from section \ref{Measures} shows that the three restricted
operators
\begin{equation}
F_{i}\text{:}=S_{i}^{\ast}\left\vert \underset{\mathcal{M}}{}\right.  ,i=0,1,2
\label{WavFracEq8}%
\end{equation}
have the following matrix representation
\begin{equation}
F_{0}=\left(
\begin{array}
[c]{ccc}%
\frac{1}{\sqrt{2}} & 0 & 0\\
0 & \frac{1}{\sqrt{2}} & 0\\
0 & 0 & 0
\end{array}
\right)  ,\;F_{1}=\left(
\begin{array}
[c]{ccc}%
0 & 0 & 0\\
0 & 0 & 1\\
0 & 0 & 0
\end{array}
\right)  ,\;F_{2}=\left(
\begin{array}
[c]{ccc}%
\frac{1}{\sqrt{2}} & 0 & 0\\
0 & \frac{-1}{\sqrt{2}} & 0\\
0 & 0 & 0
\end{array}
\right)  \text{.} \label{WavFracEq9}%
\end{equation}
A substitution of (\ref{WavFracEq8})--(\ref{WavFracEq9}) into
(\ref{WavFracEq4}) yields
\begin{align}
&  \mu_{0}\left(  \left[  \frac{i_{1}}{N}+\cdots+\frac{i_{k}}{N^{k}}%
,\frac{i_{1}}{N}+\cdots+\frac{i_{k}}{N^{k}}+\frac{1}{N^{k}}\right)  \right)
\label{WavFracEq10}\\
&  =\left\Vert F_{i_{k}}\cdots F_{i_{2}}F_{i_{1}}e_{0}\right\Vert
^{2}\nonumber\\
&  =\left\{
\begin{array}
[c]{l}%
0\text{ if one or more of the }i_{j}\text{'s is }1\\
2^{-k}\text{ otherwise.}%
\end{array}
\right. \nonumber
\end{align}
Using this, it is immediate to see that $\mu_{0}$ satisfies Hutchinson's
identity (\ref{WavFracEq7}), and therefore $\mu_{0}$ is the Hausdorff measure
$\mu^{(s)},$ $s=\log_{3}\left(  2\right)  $, restricted to $X_{3}$.

\section{A Technical Lemma\label{TechLem}}

In the proof of Lemma \ref{MeasLem2}, Theorem \ref{MeasTheo1} and Corollary
\ref{MeasCor1} above, we relied on the following lemma regarding operators in
a finite-dimensional Hilbert space. While it is analogous to the classical
Perron-Frobenius theorem, our present result makes no mention of positivity.
In fact, our matrix entries will typically be complex.

\begin{notation}
\label{TechLemNote1}If $\mathcal{M}$ is a complex Hilbert space, we denote by
$L\left(  \mathcal{M}\right)  $ the algebra of all bounded linear operators on
$\mathcal{M}$ If $\mathcal{M}$ is also finite-dimensional, we will pick
suitable matrix representations for operators $F$\emph{:}$\mathcal{M}%
\longrightarrow\mathcal{M}$ If $\mathcal{M}$ contains subspaces,
$\mathcal{M}_{i}$, $i=1,2$ such that
\[
\mathcal{M}_{1}\bot\mathcal{M}_{2}\text{ and }\mathcal{M}=\mathcal{M}%
_{1}\oplus\mathcal{M}_{2}\text{,}%
\]
then we get a block-matrix representation
\begin{equation}
F=\left(
\begin{array}
[c]{cc}%
A & B\\
C & D
\end{array}
\right)  \label{TechLemEq1}%
\end{equation}
where the entries are linear operators specified as follows.
\[
A\text{\emph{:}}\mathcal{M}_{1}\longrightarrow\mathcal{M}_{1}\text{,\quad
}B\text{\emph{:}}\mathcal{M}_{2}\longrightarrow\mathcal{M}_{1}\text{;}%
\]
and
\[
C\text{\emph{:}}\mathcal{M}_{1}\longrightarrow\mathcal{M}_{2}\text{,\quad
}D\text{\emph{:}}\mathcal{M}_{2}\longrightarrow\mathcal{M}_{2}\text{.}%
\]
If $\dim\mathcal{M}_{1}=1$, and $\mathcal{M}_{1}=\mathbb{C}w$ for some
$w\in\mathcal{M}$, then we will identify the operators $\mathcal{M}%
_{1}\longrightarrow\mathcal{M}$ with $\mathcal{M}$ via
\[
T_{\eta}\emph{:}\mathbb{C}\ni z\longrightarrow z\eta,
\]
where $\eta\in\mathcal{M}$. The adjoint operator is
\[
T_{\eta}^{\ast}x=\left\langle \eta\mid x\right\rangle w\text{, for }%
x\in\mathcal{M}.
\]

\end{notation}

\begin{lem}
\label{TechLemLem1}Let $\mathcal{M}$ be a finite-dimensional complex Hilbert
space, with $d=\dim\mathcal{M}$. Let $F\in L\left(  \mathcal{M}\right)  $, and
let $a\in\mathbb{C}$ satisfy the following four conditions\emph{:~}

\noindent\textup{(}i\textup{)} $a\in\operatorname*{spec}\left(  F\right)  ;$ ~

\noindent\textup{(}ii\textup{)} $\left\vert a\right\vert >\max\left\{
\left\vert \lambda\right\vert \mid\lambda\in\operatorname*{spec}\left(
F\right)  \backslash\left\{  a\right\}  \right\}  ;$ ~

\noindent\textup{(}iii\textup{)} the algebraic multiplicity of $a$ is one$;$
and ~

\noindent\textup{(}iv\textup{)} there is a $w\in\mathcal{M}$, $\left\Vert
w\right\Vert =1$, such that $F^{\ast}w=\bar{a}w$.

Then there is a unique $\xi\in\mathcal{M}$ such that
\begin{equation}
\left\langle w\mid\xi\right\rangle =1\text{ and }F\xi=a\xi\text{.}
\label{TechLemEq2}%
\end{equation}
Moreover,
\begin{equation}
\lim_{n\longrightarrow\infty}a^{-n}F^{n}x=\left\langle w\mid x\right\rangle
\xi\text{ for all }x\in\mathcal{M}\text{.} \label{TechLemEq3}%
\end{equation}

\end{lem}

\begin{lem}
\label{TechLemLem2}There is a constant $C$ independent of $d=\dim\mathcal{M}$
and of $x$, such that
\begin{equation}
\left\Vert a^{-n}F^{n}x-\left\langle w\mid x\right\rangle \xi\right\Vert \leq
Cn^{d-1}\max\left\{  \left\vert \frac{s}{a}\right\vert ^{n}\mid s\in
\operatorname*{spec}\left(  F\right)  \backslash\left\{  a\right\}  \right\}
\text{.} \label{TechLemEq4}%
\end{equation}

\end{lem}

\begin{pf}
$[$Lemma 5.2$]$\emph{ }Set
\[
\mathcal{M}^{\prime}\text{:}=\mathcal{M}\ominus\mathbb{C}w=\left\{
x\in\mathcal{M}\mid\left\langle w\mid x\right\rangle =0\right\}  \text{.}%
\]
Then
\begin{equation}
\mathcal{M}=\mathbb{C}w\oplus\mathcal{M}^{\prime}\text{,} \label{TechLemEq5}%
\end{equation}
and we get the resulting block-matrix representation of $F$,
\begin{equation}
F=\left(
\begin{array}
[c]{ccc}%
a & \; & 00\cdots0\\
\eta &  & G
\end{array}
\right)  \label{TechLemEq6}%
\end{equation}
where $a$ is the number in (i), the vector $\eta\in\mathcal{M}^{\prime}$, and
operator $G\in L\left(  \mathcal{M}^{\prime}\right)  $, are uniquely
determined. As a result, we get the factorization
\begin{equation}
\det\left(  \lambda-F\right)  =\left(  \lambda-a\right)  \det\left(
\lambda-G\right)  \label{TechLemEq7}%
\end{equation}
for the characteristic polynomial. Assumptions (ii) and (iii) imply
\begin{equation}
\operatorname*{spec}\left(  F\right)  \backslash\left\{  a\right\}
=\operatorname*{spec}\left(  G\right)  \text{;} \label{TechLemEq8}%
\end{equation}
and in particular, we note that $a$ is not in the spectrum of $G$. Hence the
inverse $\left(  a-G\right)  ^{-1}$ is well defined, and $\left(  a-G\right)
^{-1}\in L\left(  \mathcal{M}^{\prime}\right)  .$ We claim that the vector
\begin{equation}
\xi=w+\left(  a-G\right)  ^{-1}\eta\label{TechLemEq9}%
\end{equation}
satisfies the conditions in \textup{(}\ref{TechLemEq2}\textup{)}. First note
that $\left(  a-G\right)  ^{-1}\eta\in\mathcal{M}^{\prime}$, so $\left\langle
w\mid\xi\right\rangle =\left\langle w\mid w\right\rangle =\left\Vert
w\right\Vert ^{2}=1$. Moreover,
\[
F\xi=aw+\eta+G\left(  a-G\right)  ^{-1}\eta=aw+a\left(  a-G\right)  ^{-1}%
\eta=a\xi\text{,}%
\]
which proves the second condition in \textup{(}\ref{TechLemEq2}\textup{)}.
Uniqueness of the vector $\xi$ in \textup{(}\ref{TechLemEq2}\textup{)} follows
from \textup{(}\ref{TechLemEq8}\textup{)}. Using the matrix representation
\textup{(}\ref{TechLemEq6}\textup{)}, we get
\[
F^{2}=\left(
\begin{array}
[c]{cc}%
a^{2} & 00\cdots0\\
a\eta+G\eta & G^{2}%
\end{array}
\right)
\]
and by induction,
\begin{align}
F^{n}  &  =\left(
\begin{array}
[c]{cc}%
a^{n} & 00\cdots0\\
a^{n-1}\eta+a^{n-2}G\eta+\cdots+G^{n-1}\eta & G^{n}%
\end{array}
\right) \label{TechLemEq10}\\
&  =\left(
\begin{array}
[c]{cc}%
a^{n} & 00\cdots0\\
\left(  a^{n}-G^{n}\right)  \left(  a-G\right)  ^{-1}\eta & G^{n}%
\end{array}
\right)  \text{.}\nonumber
\end{align}
Hence, if we show that
\begin{equation}
\lim_{n\longrightarrow\infty}a^{-n}G^{n}=0\text{,} \label{TechLemEq11}%
\end{equation}
then the desired conclusion \textup{(}\ref{TechLemEq3}\textup{)} will follow.
Using the matrix form \textup{(}\ref{TechLemEq10}\textup{)}, the conclusion
\textup{(}\ref{TechLemEq3}\textup{)} reads
\begin{equation}
\lim_{n\longrightarrow\infty}a^{-n}F^{n}=\left(
\begin{array}
[c]{cc}%
1 & 00\cdots0\\
\left(  a-G\right)  ^{-1}\eta & 0
\end{array}
\right)  \text{.} \label{TechLemEq12}%
\end{equation}
In proving \textup{(}\ref{TechLemEq11}\textup{)}, we will make use of the
Jordan-form representation for $G$. Jordan's theorem applied to $G$ yields
three operators $D,V,N\in L\left(  \mathcal{M}^{\prime}\right)  $ with the
following properties:~

\noindent\textup{(}1\textup{)} $D$ is a diagonal matrix with the numbers
$\operatorname*{spec}\left(  F\right)  \backslash\left\{  a\right\}  $ down
the diagonal$;$~

\noindent\textup{(}2\textup{)} $V$ is invertible$;$ ~

\noindent\textup{(}3\textup{)} $N$ is nilpotent$\emph{:}$ If $d-1=\dim\left(
\mathcal{M}^{\prime}\right)  $ then $N^{d-1}=0;$ ~

\noindent\textup{(}4\textup{)} $\left[  N,D\right]  =ND-DN=0;$ ~

\noindent\textup{(}5\textup{)} $G=V\left(  D+N\right)  V^{-1}$. $~$

Let $x\in\mathcal{M}^{\prime}$, and let $n\geq d$. Using (2)--(5), we get
\begin{align*}
a^{-n}G^{n}x  &  =Va^{-n}\left(  D+N\right)  ^{n}V^{-1}x\\
&  =\sum_{i=0}^{d-2}\binom{n}{i}Va^{-n}D^{n-i}N^{i}V^{-1}x\text{.}%
\end{align*}
But the matrix $a^{-n}D^{n-i}$ is diagonal with entries
\[
\left\{  a^{-n}s^{n-i}\mid s\in\operatorname*{spec}\left(  F\right)
\backslash\left\{  a\right\}  \right\}  \text{, }0\leq i\leq d-1\text{.}%
\]
Using finally assumption \textup{(}ii\textup{)}, we conclude that
\begin{equation}
\lim_{n\longrightarrow\infty}\binom{n}{i}a^{-n}s^{n-i}=0\text{,}
\label{TechLemEq13}%
\end{equation}
and the proof of \textup{(}\ref{TechLemEq11}\textup{)} is completed.
\end{pf}

\begin{pf}
$[$Lemma 5.3$]$ Let the conditions be as stated in the Remark. From the
arguments in the proof of Lemma \ref{TechLemLem1}, we see that the two vectors
on the left-hand side in \textup{(}\ref{TechLemEq4}\textup{)} may be
decomposed as follows:
\begin{equation}
a^{-n}F^{n}x=\left\langle w\mid x\right\rangle w+\left(  1-a^{-n}G^{n}\right)
\left(  a-G\right)  ^{-1}\eta+a^{-n}G^{n}P_{\mathcal{M}^{\prime}}x\text{;}
\label{TechLemEq14}%
\end{equation}
and
\begin{equation}
\left\langle w\mid x\right\rangle \xi=\left\langle w\mid x\right\rangle
w+\left(  a-G\right)  ^{-1}\eta\text{.} \label{TechLemEq15}%
\end{equation}
Hence, the difference is in $\mathcal{M}^{\prime}$, and
\begin{align*}
&  \left\Vert a^{-n}F^{n}x-\left\langle w\mid x\right\rangle \xi\right\Vert \\
&  =\left\Vert a^{-n}G^{n}\left(  P_{\mathcal{M}^{\prime}}x-\left(
a-G\right)  ^{-1}\eta\right)  \right\Vert \\
&  \leq Cn^{d-1}\max\left\{  \left\vert \frac{s}{a}\right\vert ^{n}\mid
s\in\operatorname*{spec}\left(  F\right)  \backslash\left\{  a\right\}
\right\}
\end{align*}
which is the desired conclusion \textup{(}\ref{TechLemEq4}\textup{)}.
\end{pf}


\begin{thebibliography}{99}                                                                                               %


{}

\bibitem {BJMP04}L. W. Baggett, P. E. T. Jorgensen, K. D. Merrill, J. A.
Packer, An Analogue of Bratteli-Jorgensen Loop Group Actions for GMRA's, Vol.
345 of Wavelets, frames and operator theory, Contemp. Math., , American
Mathematical Society, Providence, RI, 2004, pp. 11--25.

\bibitem {BraEvaJor2000}O. Bratteli, D. E. Evans, P. E. T. Jorgensen,
Compactly Supported Wavelets and Representations of the Cuntz Relations, Appl.
Comp. Harmon. Anal., \textbf{8} (2000) 166--196.

\bibitem {BrJo02b}O. Bratteli, P. E. T. Jorgensen, Wavelets through a Looking
Glass: The World of the Spectrum, Applied and Numerical Harmonic Analysis,
Birkh\"{a}user, Boston, 2002.

\bibitem {BrJo99}O. Bratteli, P. E. T. Jorgensen, Iterated function systems
and permutation representations of the Cuntz algebra, Mem. Amer. Math. Soc.
139 (1999) no. 663.

\bibitem {CMW92}R. Coifman, Y. Meyer, V. Wickerhauser, Size properties of
wavelet-packets, in: M. Ruskai, G. Beylkin, R. Coifman, I. Daubechies, S.
Mallat, Y. Meyer, L. Raphael (Eds.), Wavelets and Their Applications, Jones
and Bartlett, Boston, 1992, pp. 453--470.

\bibitem {Cu77}J. Cuntz, Simple $C^{\ast}$-algebras generated by isometries,
Comm. Math. Phys. 57 (1977) 173--185.

\bibitem {Dau92}I. Daubechies, Ten lectures on wavelets, CBMS-NSF Regional
Conference Series in Applied Mathematics, Philadelphia: Society for Industrial
and Applied Mathematics,Vol. 61, 1992.

\bibitem {DuJo05}D. E. Dutkay; P. E. T. Jorgensen, Wavelet constructions in
non-linear dynamics. Electron. Res. Announc. Amer. Math. Soc. 11. (2005),
21--33 (electronic).

\bibitem {DutJor04}D. Dutkay, P. Jorgensen, Wavelets on Fractals, preprint,
University of Iowa, 2004, arXiv:math.CA/0305443, accepted: Rev. Mat. Iberoamericana.

\bibitem {Falconer85}K. J. Falconer, The Geometry of Fractal Sets, Cambridge,
Tracts on Mathematics, \textbf{85}.

\bibitem {Hut81}J. E. Hutchinson, Fractals and self-similarity, Indiana Univ.
Math. J. \textbf{30} (1981) 713--747.

\bibitem {JoKr03}P. E. T. Jorgensen; D. W. Kribs, Wavelet representations and
Fock space on positive matrices. J. Funct. Anal. 197 (2003), no. 2, 526--559.

\bibitem {Jor01}P.E.T. Jorgensen, Minimality of the Data in Wavelet Filters,
Advances in Mathematics \textbf{159}, 2001, pp. 143--228.

\bibitem {Jor04}P.E.T. Jorgensen, Measures in wavelet decompositions, Adv.
Appl. Math., Vol. 23 (2005), 561--590.

\bibitem {Wic93}M. V. Wickerhauser, Best-adapted wavelet packet bases, in: I.
Daubechies (Ed.), Different Perspectives on Wavelets (San Antonio, TX, 1993),
Vol. 47 of Proc. Sympos. Appl. Math., American Mathematical Society,
Providence, 1993, pp. 155--171.

\bibitem {Wic94}M. V. Wickerhauser, Adapted Wavelet Analysis from Theory to
Software, IEEE Press, New York, A.K. Peters, Wellesley, MA, 1994.
\end{thebibliography}
\end{document}